\documentclass[10pt]{article}   	

\usepackage{amsmath,amsthm, amssymb,amsfonts,amscd,mathabx}
\usepackage[title]{appendix}
\usepackage[colorlinks=true,linkcolor=blue,citecolor=red]{hyperref}
\usepackage{mathrsfs}                             
\usepackage[pdftex]{graphicx}                 
\usepackage{float}
\usepackage{mathtools}
\usepackage{tikz}
\usetikzlibrary{matrix}
 \let\overbrace\LaTeXoverbrace
\numberwithin{equation}{section}
\newtheorem{theorem}{Theorem}[section]
\newtheorem{lemma}[theorem]{Lemma}
\newtheorem{proposition}[theorem]{Proposition}
\newtheorem{definition}[theorem]{Definition}
\newtheorem{corollary}[theorem]{Corollary}
\newtheorem{hypothesis}[theorem]{Hypothesis}
\newtheorem{remark}[theorem]{Remark}

\newenvironment{Proof}[1][.]%
{\begin{trivlist}\item[]\textbf{Proof#1 }}%
{\qed\end{trivlist}}

 {\begin{trivlist}\item[]\textbf{Acknowledgments }}{\end{trivlist}}

\newcommand{\bt}{\begin{theorem}}
\newcommand{\et}{\end{theorem}}
\newcommand{\bl}{\begin{lemma}}
\newcommand{\el}{\end{lemma}}
\newcommand{\bp}{\begin{proposition}}
\newcommand{\ep}{\end{proposition}}
\newcommand{\bd}{\begin{definition}}
\newcommand{\ed}{\end{definition}}
\newcommand{\bc}{\begin{corollary}}
\newcommand{\ec}{\end{corollary}}
\newcommand{\br}{\begin{remark}}
\newcommand{\er}{\end{remark}}
\newcommand{\bh}{\begin{hypothesis}}
\newcommand{\eh}{\end{hypothesis}}
\newcommand{\be}{\begin{enumerate}}
\newcommand{\ee}{\end{enumerate}}
\newcommand{\beq}{\begin{equation}}
\newcommand{\eeq}{\end{equation}}
\newcommand{\beqs}{\begin{equation*}}
\newcommand{\eeqs}{\end{equation*}}
\newcommand{\bpf}{\begin{Proof}}
\newcommand{\epf}{\end{Proof}}
\newcommand{\bld}{\begin{aligned}}
\newcommand{\eld}{\end{aligned}}
\newcommand{\bhp}{\begin{hypothesis}}
\newcommand{\ehp}{\end{hypothesis}}
\newcommand{\bcs}{\begin{cases}}
\newcommand{\ecs}{\end{cases}}
\newcommand{\R}{\mathbb{R}}

\newcommand{\N}{\mathbb{N}}
\newcommand{\Z}{\mathbb{Z}}

\newcommand{\T}{\mathbb{T}}

\newcommand{\rmnum}[1]{\romannumeral #1}
\newcommand{\Rmnum}[1]{\uppercase\expandafter{\romannumeral #1\relax}}

\newcommand{\caO}{\mathcal{O}}
\newcommand{\rmO}{\mathrm{O}}

\newcommand{\rmd}{\mathrm{d}}
\newcommand{\rme}{\mathrm{e}}
\newcommand{\rmi}{\mathrm{i}}

\newcommand{\id}{\mathrm{\,Id}}

\newcommand{\qprime}{{\prime\prime\prime\prime}}
\newcommand{\tprime}{{\prime\prime\prime}}
\newcommand{\dprime}{{\prime\prime}}

\renewcommand{\Re}{\mathrm{Re}}
\renewcommand{\Im}{\mathrm{Im}}

\renewcommand{\leq}{\leqslant}
\renewcommand{\geq}{\geqslant}

\def\ker{\mathop{\mathrm{\,Ker}\,}}

\def\rg{\mathop{\mathrm{\,Rg}\,}}

\def\veps{\vepsilon}
\def\veps{\varepsilon}
\def\vphi{\varphi}

\def\Xint#1{\mathchoice
   {\XXint\displaystyle\textstyle{#1}}%
   {\XXint\textstyle\scriptstyle{#1}}%
   {\XXint\scriptstyle\scriptscriptstyle{#1}}%
   {\XXint\scriptscriptstyle\scriptscriptstyle{#1}}%
   \!\int}
\def\XXint#1#2#3{{\setbox0=\hbox{$#1{#2#3}{\int}$}
     \vcenter{\hbox{$#2#3$}}\kern-.5\wd0}}

\def\dashint{\Xint-}

\def\caD{\mathcal{D}}

\newcommand{\caL}{\mathcal{L}}

\newcommand{\caN}{\mathcal{N}}
\newcommand{\caG}{\mathcal{G}} 
\newcommand{\caR}{\mathcal{R}}
\newcommand{\caB}{\mathcal{B}}
\newcommand{\caM}{\mathcal{M}}
\newcommand{\caH}{\mathcal{H}}
\newcommand{\caP}{\mathcal{P}}

\newcommand{\unu}{\underline{\nu}}

\newcommand{\ux}{\underline{x}}
\newcommand{\uy}{\underline{y}}
\newcommand{\ua}{\underline{a}}

\newcommand{\LambdaO}{\Lambda_\mathrm{eff}}

\makeatletter\@addtoreset{figure}{section}\makeatother

\newfam\bifam
\font\tenbi=cmmib10 scaled \magstep1 \font\sevenbi=cmmib10 at 11pt
\font\fivebi=cmmib10 at 6pt \textfont\bifam = \tenbi
\scriptfont\bifam = \sevenbi \scriptscriptfont\bifam= \fivebi

\begin{document}

\begin{center}
{\fontsize{15}{15}\fontfamily{cmr}\fontseries{b}\selectfont{Pinning and dipole asymptotics of locally deformed striped phases
}}\\[0.2in]
Arnd Scheel$\,^{\dagger,1}$ and Qiliang Wu$\,^{\ddagger,}$\footnote{The authors acknowledge partial support by the National Science Foundation through grants  NSF DMS-2205663 (AS) and  NSF DMS-1815079 (QW). } \\
\textit{\footnotesize $\,^\dagger$University of Minnesota, School of Mathematics,   206 Church St. S.E., Minneapolis, MN 55455, USA\\
$\,^\ddagger$Ohio University, Department of Mathematics,  Morton Hall 537, 1 Ohio University, Athens, OH 45701, USA}\\[3mm]
{\small \today}
\end{center}

\begin{abstract}
\noindent 
We investigate the effect of spatial inhomogeneity on perfectly periodic, self-organized striped patterns in spatially extended systems. We demonstrate that inhomogeneities select a specific translate of the striped patterns and induce algebraically decaying, dipole-type farfield deformations. Phase shifts and leading order terms are determined by effective moments of the spatial inhomogeneity. Farfield decay is proportional to the derivatives of the Green's function of an effective Laplacian. Technically, we use mode filters and conjugacies to an effective Laplacian to establish Fredholm properties of the linearization in Kondratiev spaces. Spatial localization in a contraction argument is gained through the use of an explicit deformation ansatz and a subtle cancellation in Bloch wave space. 
\end{abstract}

\section{Introduction}

Striped patterns are one of the simplest self-organized structures in spatially extended systems, with observations and analysis including sand \cite{sand} and icicle ripples \cite{icicle},  convection roll  \cite{bodenschatz} and precipitation patterns \cite{thomas},  bacterial colony growth \cite{eshel} or the formation of presomites in early development \cite{digit},   ion-beam milling \cite{bradley}, dip-coating \cite{dipstripe}, lamellar crystal growth \cite{double,zigzageutectic}, wrinkling patterns \cite{aharoni2017smectic}, and water jet cutting \cite{friedrich}. Although the physical mechanisms leading to the formation of stripes vary, one can often find universal features in the dynamics associated with systems that exhibit striped patterns. Emergence of nonlinear striped patterns is often associated with an instability that can be studied as a simple symmetry breaking bifurcation in a system with translation symmetry \cite{turing1952chemical,chossatlauterbach,chossatiooss}. The assumptions here are that the system possesses an idealized translation invariance and that the  observed pattern possesses an idealized spatial periodicity. Actually observed patterns rarely exhibit these perfect features: systems are not translation invariant due to boundaries or spatial inhomogeneities  
\cite{goh_scheel_review,scheelweinburd,scheelmorrissey,JSW,JaramilloScheel_15}, and patterns typically exhibit many self-organized defects, such as dislocations, disclinations, or grain boundaries
\cite{zanv_21,ercolani_etal_2001,sw_14,LloydScheel,pismen_book}. Nevertheless, there have been many successful attempts at identifying universal features of such imperfect striped patterns, including descriptions of defects, boundary conditions, and the effect of inhomogeneities. 

A prototypical model for the dynamics of stripes is the Swift-Hohenberg equation, originally introduced to mimic convection roll dynamics in B\'enard convection \cite{sh} but used across the sciences, early in Turing's notes on morphogenesis \cite{dawes}, and in applications from   plant phyllotaxis \cite{phyllo} to  nonlinear optics \cite{shoptics},
\begin{equation}\label{e:sh0}
u_t=-(\Delta+1)^2u+\mu u - u^3,\qquad x\in\R^n,\ u\in\R,\ t\geq0,
\end{equation}
where the linear equation, at $\mu\gtrsim 0$,  amplifies Fourier modes with wavenumber $|\underline{k}|\sim 1$ but dampens all other  modes. As a result, one observes solutions that locally in space resemble a sinusoidal pattern $u\sim A\sin(\underline{k}\cdot \underline{x}-x_0)$ for some amplitude $A$, phase shift $x_0$, and wave vector $\underline{k}\in\R^n$ with $|\underline{k}|\sim 1$. 

One can proof existence of periodic solutions $u_*(\underline{k}\cdot \ux;|\underline{k}|)$ for the nonlinear equation for parameter values $\mu\gtrsim 0$ using bifurcation theory. One can also study temporal stability of such solutions for perturbed initial conditions, starting with the linearization. The most informative results on stability are concerned with a setting where $\ux\in\R^n$. Solutions are linearly stable within a strict subset of the existence region in the $(\mu,|\underline{k}|)$-plane, with boundaries including zigzag and Eckhaus instabilities in this simple setting; see for instance \cite{crosshohenberg,mielke_1997}. In the interior and to some extent on the boundary of the stability region, one can proof decay of (non-periodic) perturbations to a striped pattern: starting with initial conditions of the form $u_*+w$ at time $t=0$, where $u_*$ is a perfect striped patterns and $w$ is small in a suitable norm, one can show that the solution $u(t,\ux)$ converges to $u_*(\ux)$ for $t\to\infty$, at least for perturbations that exhibit some degree of spatial localization; see for instance \cite{gallayscheel_2011,gswz18,johnsonzumbrun_2011,SW_15,schneider_1996,uecker_1999} and, for less localized perturbations, \cite{derijk22,jnrz13,SSSU_12}.  The key insight, that we also exploit in the present work, is that linear dynamics near striped patterns are diffusive: a homogenized equation yields an effective nonlinear diffusion equation for the wave vector, known as the Cross-Newell equation \cite{crossnewell}. Diffusion exhibits temporal decay provided nonlinearities are irrelevant, that is, of sufficiently high power, or involving derivatives that lead to cancellations \cite{bk94}.  Such cancellations also occur in the full, non-homogenized, Swift-Hohenberg equation, only in a rather hidden form \cite{schneider_1996}. The leading-order description via a diffusion equation shares many commonalities with homogenization in overdamped elasticity, where in fact many similar questions as to the effect of boundaries or the presence of defects are of interest \cite{ortner,luskinortner}. 

Temporal decay in a diffusion equation depends crucially on localization of initial conditions --- or, more generally speaking, of perturbations in space-time. Thinking of perturbations to the Swift-Hohenberg equation as given through an additional term $g(t,\ux)$ on the right-hand side of \eqref{e:sh0}, one would typically require strong temporal decay and possibly spatial decay. 

Our interest here is in the response of striped patterns to perturbations that do not decay in time, particularly perturbations of the form $g(t,\ux)\equiv g(\ux)$. One then expects the striped pattern to relax to a deformed pattern, much as a local stress induces a farfield deformation in an elastic medium. We wish to establish the existence of such a deformed striped pattern and describe the deformation at leading order. In a stationary diffusion (or Poisson) equation such a description is rather straightforward, and far-field expansions are easily obtained from moments of the inhomogeneity via multipole expansions \cite{jackson_classical_1999}. In nonlinear equations, such expansions are less obvious, although recent results establish such expansions in lattice-models for nonlinear elasticity \cite{ortner}.

Our main result provides such a leading-order expansion of deformed patterns. As a main difficulty, compared with stability results, our perturbation is less localized and some techniques available for temporal stability, such as temporal iteration or semigroup methods,  are not available. Comparing with results on scalar or lattice equations, a key difficulty here is the less explicit nature of the linearization and the necessity to detect crucial cancellations in the nonlinearity. In some ways, the most difficult case, with weakest diffusion and most nonlinear relevant terms is the case of one space dimension, which was treated in \cite{JSW}. The one-dimensional setup however does simplify the analysis  in several crucial ways, allowing for instance for the application of spatial dynamics methods. More technically, factorization of Fourier symbols is significantly easier in one space-dimension and asymptotics at spatial infinity are exponential, while, in dimension greater than 1, asymptotics are algebraic related to the fact that the dimension of the space of harmonic polynomials is infinite in dimensions greater than 1. 

Technically, we introduce a novel and direct functional analytic approach which relies on fine Fredholm mapping properties of the linear operator and bordering of the nonlinear equation with an explicit farfield ansatz. We believe that many of the techniques developed here will prove useful in the analysis of defects in crystalline media, in the context of analysis and bifurcations, as well as for computational approaches.

\paragraph{Setup and main results.}

We consider the Swift-Hohenberg equation with spatially localized impurity,
\beq\label{e:imSHE}
u_t=-(\Delta+1)^2u+\mu u - u^3 + \veps g(\ux), \qquad \ux=(x_1,\ldots,x_n)\in\R^n,n=2,3,
\eeq
where $\Delta=\nabla_{\ux}\cdot\nabla_{\ux}=\sum_{j=1}^n\partial_{x_jx_j}$ and the parameter $\mu>0$ is fixed.
Our starting point are particular solutions to \eqref{e:imSHE} at $\veps=0$, a family of striped patterns $u(t,x)=u_*(k x_1;k)=u_*(k(x_1+\frac{2\pi}{k});k)$ for wavenumbers $k\sim 1$. Striped patterns are found as solutions to the ODE boundary-value problem
\beq\label{e:k}
-\left(k^2\frac{\rmd^2 }{\rmd \xi^2}+1\right)^2u+\mu u - u^3=0, \qquad u(\xi;k)=u(\xi+2\pi;k).
\eeq
We state our results in the following for perturbations of the specific striped pattern with wavenumber $k=1$, although our results can easily be adapted to other wavenumbers. 

\begin{hypothesis}[Existence and Diffusive Stability of Stripes]\label{h:1}
We assume that there exists a stripe solution $u_*(\xi;k_*)$,
that is, $u_*(\xi;k_*)$ solves \eqref{e:k} with $k_*=1$ and $\mu>0$, fixed. We assume that the striped solution is linearly diffusively stable, that is, the effective diffusivity constants $d_\|,d_\perp$ obtained from the linearization at $u_*$ after Fourier-Bloch transform are positive; see Lemma \ref{l:1.2} for expressions for $d_\|,d_\perp$.
\end{hypothesis}

Hypothesis \ref{h:1} holds for $\mu\gtrsim 0$ \cite{mielke_1997}. Isotropy of the equation guarantees that we may focus on stripes depending on $x_1$, only, but generalizations to anisotropic settings are straightforward. To emphasize the special role of the orientation imposed by the stripe, we write
\[
\ux=(x_1, x_2, \cdots, x_n)=(x_1, \ux_\perp), \quad \ux_\perp=(x_2, \cdots, x_n).
\]
Hypothesis \ref{h:1} guarantees that the stripe solution can be extended to a smooth family of solutions $u_*(\xi;k)$, $k\sim 1$, $\xi=kx$; see for instance \cite{JSW}. Define
\[
\begin{matrix}
u_*^\prime:=\frac{\partial u_*(\xi;k)}{\partial \xi}\mid_{k=k_*}, & u_{*,k}:=\frac{\partial u_*(\xi;k)}{\partial k}\mid_{k=k_*}, & u_{*,k}^\prime:=\frac{\partial^2 u_*(\xi;k)}{\partial \xi\partial k}\mid_{k=k_*},
\end{matrix}
\]
with similar notation for higher order derivatives, for instance  $u_*^\dprime:=\partial_\xi^2 u_*\mid_{k=k_*}$. We will omit the parameter $k$ in the notation for $u_*$ when $k=1$, that is, $u_*(\xi):=u_*(\xi;1)$.

Clearly, $u_*(k(x_1+x_1^0);k)$ is a solution for any $x_1^0$ as well, providing us with a neutral mode in the equation. Our second hypothesis is concerned with the projection of the inhomogeneity onto this neutral mode. Define therefore 
\beq\label{e:melnikov}
    M(x_1^0):=\int_{\R^n} u_*'(x_1+x_1^0)g(\ux)d \ux. 
\eeq
Note that  $\int_0^{\frac{2\pi}{k}} M(x_1^0)\rmd x_1^0=0$, and $M(x_1^0)=M(x_1^0+\frac{2\pi}{k})$, 
so that, for fixed $k$, there exists an (in fact at least two) $x_1^0=:x_1^*$ where $M(x_1^*)=0$.
We shall assume that at least one of those zeros is non-degenerate.
\begin{hypothesis}[Non-degenerate pinning]\label{h:2}
Fix $k=k_*$ from Hypothesis \ref{h:1}. We assume that for some $x_1^*\in [0,2\pi/k)$, we have 
\[
M(x_1^*)=0,\qquad 
M^\prime(x_1^*)\neq 0.
\]
\end{hypothesis}
Our main result will show that perturbed stripe patterns are ``pinned at'' such positions $x_1^*$, motivating our terminology of pinning. The results in \cite{JaramilloScheel_15} for the Ginzburg-Landau equation suggest that stationary patterns ``centered'' at different locations are sustained only with logarithmically growing farfield corrections. 

In order to state our main result, we introduce algebraically weighted Sobolev spaces. Let $L^p(\R^n)$ be the usual Lebesgue spaces and  $W^{\ell,p}(\R^n)$ the Sobolev spaces with $\ell^\text{\footnotesize{th}}$ derivative in $L^p(\R^n)$, and $H^\ell(\R^n)=W^{\ell,2}(\R^n)$, where we usually drop the argument $\R^n$. Define $\langle \ux\rangle=(1+|\ux|^2)^{1/2}$ and the weighted spaces $L^p_\gamma$, $W^{\ell,p}_\gamma$ as the closure of $C^\infty_0$ in the norms
\[
\|w\|_{L^p_\gamma}:=\|w(\cdot) \langle \cdot\rangle^\gamma\|_{L^p},\qquad 
\|w\|_{W^{\ell,p}_\gamma}:=\|w(\cdot) \langle \cdot\rangle^\gamma\|_{W^{\ell,p}},
\]
and again $H^\ell_\gamma:=W^{\ell,2}_\gamma$.
In order to describe farfield deformations, we need the 
Green's function to the anisotropic Laplacian $\Delta_\mathrm{eff}=d_\|\partial_{x_1}^2+d_\perp\Delta_{x_\perp}$, which we denote by $G(\ux)$; see \eqref{e:deltaeff}. We also introduce a smoothed characteristic function $\chi\in C^\infty_0(\R^n)$,
\begin{equation}\label{e:chi}
\chi(\ux)=1 \text{ for }|\ux|<1,\qquad  \chi(\ux)=0 \text{ for } |\ux|\geq 2,\qquad \chi(\ux)\in[0,1] \text{ for all } \ux\in\R^n,
\end{equation}
as well as the scaled versions
\beq\label{e:rechi}
\chi_r(\ux):=\chi(\ux/r), \quad \text{for any }r>0.
\eeq

\begin{theorem}[Deformation of stripes]\label{t:main}
Fix the space dimension $n\in \{2,3\}$. Assume Hypothesis \ref{h:1} on existence and linear diffusive stability and Hypothesis \ref{h:2} on non-degenerate pinning for some $x_1^*$. Fix $2+\frac{n}{2}>\gamma>1+\frac{n}{2}$ and assume that the inhomogeneity is smooth and spatially localized, $g\in H^2_\gamma$. Then there exists $\veps_0>0$ such that for all $0<\veps<\veps_0$ there exists a stationary  solution 
\[
u(\ux)=u_*(x_1+a_0(\veps)+\Theta(\ux;\ua(\veps));\kappa(\ux;\ua(\veps))+h(\ux;\veps)u_*'(x_1+a_0(\veps))+w(\ux;\veps),
\]
to \eqref{e:imSHE} with associated phase and wavenumber deformation 
\[
\Theta(\ux;\ua(\veps))=\ua(\veps)\cdot\nabla_{\ux} \big((1-\chi(\ux))G(\ux)\big),\qquad 
\kappa(\ux;\ua(\veps))=|1+\nabla_{\ux} \Theta(\ux;\ua(\veps))|,
\]
and strongly localized  corrections 
$h\in L^2_{\gamma-2}$ and $w\in L^2_{\gamma-1}$. Moreover,  $a_0$, $\ua=(a_1,\cdots,a_n)$,  $h\in L^2_{\gamma-2}$, and $w\in L^2_{\gamma-1}$ are smooth  in $\veps$, with leading order expansions with $h(\ux;\veps)=\caO(\veps)$ and $ w(\ux;\veps)=\caO(\veps)$, and
\[
\begin{aligned}
&a_0(\veps)=x_1^*+\caO(\veps), \qquad a_j(\veps)=\veps \,
\frac{
\int_{\R^n} H_{1,j}(\ux)g(x_1-x_0^*,\ux_\perp)\rmd \ux
}
{(d_{||}d_{\perp}^{n-1})^{1/2 }
\dashint_0^{2\pi} (u_*^\prime)^2
}+\caO(\veps^2)
\text{ for }1\leq j\leq n,
\end{aligned}
\]
where we used the notation for the average of a periodic function $\dashint v=\frac{1}{2\pi}\int_0^{2\pi}v(x_1)\rmd x_1$, and $H_{1,j}$ are the pseudo-harmonic polynomials of degree 1 
\begin{align}
  H_{1,1}(\ux)=x_1u_*^\prime(x_1)+u_{*,k}(x_1), \qquad   H_{1,\ell}(\ux)=x_\ell u_*^\prime(x_1), \ 2\leq \ell\leq n.\label{e:php0}
\end{align}
\end{theorem}
We note that the correction is localized. In fact, we show that in the far field, $|\ux|\gg 1$,  $1\gg|\Theta(\ux)|\gg |h(\ux)|\gg |w(\ux)|$ so that the  correction induced by the phase modulation $\Theta$ is indeed the leading-order description of the deformation. The phase correction in fact decays as $1/|\ux|^{n-1}$, while $L^2_{\gamma-1}$ corresponds to $|\ux|^\alpha$ with {$\alpha<-(\gamma-1+\frac{n}{2})<-n$} and  $L^2_{\gamma-2}$ induces decay $|x|^{\alpha+1}$,  for homogeneous decay. In fact, we also control derivatives of $h$ and $w$ which enforce this type of  decay in a pointwise sense. 
\begin{remark}[Effective dipole expansion]
    If one, somewhat artificially, sets $u_*(\xi;k)=\xi$, then the effective harmonic polynomials are in fact the linear harmonic polynomials and the farfield deformation $\Theta$ is simply the dipole field generated by the moments of $g$. This simple scenario occurs when studying a phase-diffusion or Cross-Newell approximation to the dynamics of periodic patterns, which effectively averages over small scales. As a result, in this case we do not find nondegenerate pinning as $M(x_1^0)$ is constant: the averaging destroyed the ``small-scale'' structure of the pattern and its interaction with a localized inhomogeneity.
\end{remark}
\paragraph{Outline.} The remainder of this paper is organized as follows. We review Fredholm properties of the Laplacian in weighted spaces in Section \ref{s:2} before establishing corresponding mapping properties of the linearization at a striped pattern in Section \ref{s:3}. Section \ref{s:4} introduces a phase-modulation ansatz based on Green's function expansions of the phase of the pattern and sets up a nonlinear fixed point equation. Our main result is a consequence of well-posedness of this fixed point equation,  exploiting the phase modulation, mapping properties of the linearization, and a cancellation for leading-order terms. We conclude with a discussion.

\section{Laplacian, Fredholm properties, and Kondratiev spaces}\label{s:2}

The perturbation result in our main theorem relies on a detailed understanding of the linearization at a striped pattern. The linearization is a well-defined elliptic operator on say $L^2_\gamma$ with domain $H^4_\gamma$, but it is not invertible for any $\gamma$, in fact it is not a Fredholm operator. This is in analogy with the Laplacian, which is not Fredholm but rather possesses essential spectrum at the origin. The remedy we pursue here is to work on algebraically localized spaces with an enlarged domain, contained in $L^2_{\gamma-2}$. Key to our main fixed point argument is a fine characterization of this domain and its localization properties. Before embarking on this analysis in Section \ref{s:3}, we first recall the relevant result in the case of the Laplacian. 

\subsection{Fredholm properties of $\Delta$}

We start by formally considering the kernel of the Laplacian on spaces of polynomials in $\R^n$. We write  $\caH_j(\Delta)$ for the set of  harmonic polynomials of degree $j$ and $\caH_{\leq m}(\Delta)=\oplus_{j=0}^m \caH_j(\Delta)$. For instance, when $n=2$, $\caH_m(\Delta)=\mathrm{span}\,\{\Re(x_1+i x_2)^m,\Im(x_1+i x_2)^m\}$. 
Formally exploiting self-adjointness, we expect the range and the orthocomplement of the kernel of the Laplacian
to be given by 
\[
  \caH_{\leq m}^{\perp_\gamma}(\Delta):=\left\{ f\in L^2_\gamma(\R^n) \,\middle\vert\, \int_{\R^n}fH\rmd x=0, \text{ for all } H\in \caH_{\leq m}(\Delta)\right\},
\]
with appropriate values of $m=m(\gamma)$.
Next, we encode loss of localization from inverting the Laplacian in spaces with algebraic weights, where algebraic localization is restored through differentiation. We therefore introduce the Kondratiev space $M^{k, p}_\gamma(\R^n)$ as the completion of $C_0^\infty(\R^n)$ in the norm 
\[
\|f\|_{M^{k, p}_\gamma}:=\left(\sum\limits_{|\underline{\alpha}|\leq k}\|\langle \ux\rangle^{|\underline{\alpha}|+\gamma}D^{\underline{\alpha} }f\|^p_{L^p(\R^n)}\right)^{1/p},
\]
where $\langle \ux\rangle=(1+|\ux|^2)^{1/2}$, $1<p<\infty$, $\gamma\in\R$ and $k\in\N$. For $k=0$, these are classical weighted $L^p$-spaces and thus denoted as $L^p_\gamma$. 

\begin{theorem}[Fredholm properties of $\Delta$]\cite{McOwen}
\label{t:McOwen}
Fix $n\geq 2$, $1<p<\infty$. 
Then 
\beq\label{e:lap}
\Delta: M^{2,p}_{\gamma-2}(\R^n) \longrightarrow L^p_{\gamma}(\R^n)
\eeq
is a Fredholm operator for any $\gamma$ except when $\displaystyle \gamma+n/p\in \Z\backslash (2,n)$. More specifically, 
\begin{enumerate}
\item for $\gamma+n/p\in(2, n)$, \eqref{e:lap} is an isomorphism;
\item for $\gamma+n/p\in (n+m, n+m+1)$, for some $m\in\{0,1,2,\ldots\}$, \eqref{e:lap} is an injection with closed range 
\[
R_\gamma=\caH_{\leq m}^{\perp_\gamma}(\Delta);
\]
\item for $\gamma+n/p\in(1-m,2-m)$ for some $m\in \{0,1,2,\ldots\}$, \eqref{e:lap} is a surjection with kernel
\[
N_\gamma=\caH_{\leq m}(\Delta).
\]
\end{enumerate}
On the other hand, if $\gamma+n/p=2-m$ or $\gamma+n/p=n+m$ for some $m\in \{0,1,2,\ldots\}$, then \eqref{e:lap} does not have a closed range.
\end{theorem}
\br
Throughout this subsection, we restrict ourselves to the $n$-dimensional Laplacian with $2\leq n$. For the Fredholm properties of the 1-d Laplacian, see for instance \cite{JSW}.
\er
We will make use of an equivalent of the Laplacian, separating regularization from localization.

\bp
The closed operator 
\[
\Delta(1-\Delta)^{-1}: L^p_{\gamma-2}(\R^n)\longrightarrow L^p_\gamma(\R^n)
\]
has the same Fredholm properties as the Laplacian stated in Theorem \ref{t:McOwen}, that is, Fredholm indices and cokernels coincide and kernels are conjugate with the isomorphism $(1-\Delta)^{-1}$.
\ep
The proposition is an immediate consequence of the fact that 
\[
1-\Delta: W^{2,p}_{\gamma}(\R^n)\longrightarrow L^p_\gamma(\R^n)
\]
is bounded and invertible for any $1<p<\infty$ and $\gamma\in\R$, and commutes with derivatives. We omit details. 

We conclude with some simple consequences. 
\br
\begin{enumerate}\label{r:owendis}
    \item Scalings $\ux\mapsto D\ux$ with $D=\mathrm{diag}(d_j)$, $d_j>0$, induce anisotropic versions of the Laplacian $D\Delta D^{-1}$, which possess the same Fredholm properties, that is, they have the same Fredholm index and respective kernels and cokernels obtained via scaling from the harmonic polynomials of the Laplacian. 
    \item Spatially discrete derivatives improve localization as well, a fact that will be useful in the context of operators with spatially periodic coefficients. If we write  $\delta_j$ for the discrete derivative, $(\delta_j u)(x)=u(x+2\pi e_j)-u(x)$, then 
    \[
        \|D^\beta \delta^\alpha \Delta^{-1} f\|_{L^2_{\gamma-2+|\beta|+|\underline{\alpha}|}}\leq \|f\|_{L^2_\gamma},
    \]
    for all multi indices $\underline{\alpha},\beta$ with $|\underline{\alpha}|+|\beta|\leq 2$ and $f$ in the range of $\Delta$. To see this, simply notice that $\delta_j [\partial_{x_j} (1-\partial_{x_j})^{-1}]^{-1}$ is bounded on $L^2_\gamma$, which in turn follows immediately from inspecting the Fourier symbol.
\end{enumerate}
\er

\section{Mapping properties of the linearization at stripes}\label{s:3}
Studying perturbations of $u_*(x_1)$ as a solution to the stationary Swift-Hohenberg equation, one is first led to examining properties of the linearized operator 
\begin{equation}
    \caL_*:=-(\Delta +1)^2 u + \mu u -3u_*^2 u.
\end{equation}
Clearly, $\caL_*$ is closed, elliptic, with domain of definition $H^4_\gamma$ when considered on $L^2_\gamma$,  in fact self-adjoint for $\gamma=0$. Our goal here is to obtain precise characterizations of mapping properties of $\caL_*$. We therefore start in Section \ref{s:3.1} by  introducing a Bloch wave description and finding pseudo-harmonic polynomials representing kernel and cokernel in suitable spaces. We then derive a first result showing that the inverse is bounded when giving up 2 degrees of localization, and appropriately restricting to the closed range in Section \ref{s:3.2}. The result hinges on mode filters and a conjugacy to the Laplacian. Unfortunately, the loss of localization cannot be compensated for in a nonlinear argument through gain of localization in the nonlinearity. We therefore develop a  refined description that relies on mode filters that we introduce  in Section \ref{s:3.3}, together with a conjugacy that separates a bounded invertible part of the operator from close-to-neutral modes where the action is conjugate to the Laplacian. Key to our nonlinear argument is a representation of close-to-neutral modes as envelopes of $u_*'$. We conclude in Section \ref{s:3.4} with the main linear result, mapping properties of the linearization $\caL_*$ between neutral and far-from-neutral modes, Proposition \ref{p:II}.

\subsection{Bloch wave decomposition and pseudo-harmonic polynomials of $\caL_*$}\label{s:3.1}
The linearized operator 
\[
\begin{matrix}
\caL_*:&H^4(\R^n) & \longrightarrow & L^2(\R^n) \\
& v(\ux) & \longmapsto &-(\Delta+1)^2v+\mu v-3u_*^2 v
\end{matrix}
\]
is conjugate to the direct sum of the family of Bloch-Fourier operators,
\[
\begin{matrix}
\widehat{\caL_*}(\unu):&H^4(\T_{2\pi}) & \longrightarrow & L^2(\T_{2\pi}) \\
& v(\xi) & \longmapsto &-((\partial_{\xi}+\rmi\nu_1)^2-|\unu_h|^2+1)^2v+\mu v-3u_*^2 v,
\end{matrix}
\]
where $\unu:=(\nu_1, \unu_h)\in\Omega$ with $\Omega=[-\frac{1}{2},\frac{1}{2}]\times\R^{n-1}$. 
More specifically, the Bloch-Fourier transform
\beq\label{e:BFT}
\begin{aligned}
\caB: \quad&L^2(\R^n)&\quad\longmapsto\quad &\quad L^2(\Omega, L^2(\T_{2\pi}))\\
& \quad u(\ux) &\quad \longrightarrow \quad&\quad\widetilde{u}(\unu; \xi):=\sum_{k\in\Z} \widehat{u}(k+\nu_1, \unu_h) \rme^{\rmi k \xi},
\end{aligned}
\eeq
where $\widehat{u}(\unu)=\frac{1}{(2\pi)^n}\int_{\R^n}u(\ux)\rme^{-\rmi \unu\cdot\ux}\rmd \ux$ is the Fourier transform of $u$, diagonalizes $\caL_*$, that is, 
\[
\big(\caB\circ\caL_*\circ\caB^{-1} \widetilde{u}\big)(\unu;\xi)=\widehat{\caL_*}(\unu)\widetilde{u}(\unu;\xi), \quad \forall (\unu,\xi)\in\Omega\times \T_{2\pi},
\]
denoted as $\caB\circ\caL_*\circ\caB^{-1} =\int_{\Omega}\widehat{\caL_*}(\unu)\rmd \unu$.
We have the following lemma.
\bl\label{l:1.2}
Assume that the eigenvalue $\lambda=0$ of $\widehat{\caL_*}(0)$ is algebraically simple with eigenspace spanned by $u_*'$. Then the family of Bloch-Fourier operators has the following spectral properties.
\begin{itemize}
\item[(\rmnum{1})] $\sigma(L_*)=\bigcup\limits_{\unu\in\Omega}\sigma\left(\widehat{\caL_*}(\unu)\right)$.
\item[(\rmnum{2})] $\widehat{\caL_*}(\unu)$ is bounded and invertible for $\unu\neq{ 0}$.
\item[(\rmnum{3})] There exist locally analytic continuations of the eigenvalue $\lambda(\unu)$ with $\lambda(0)=0$ and  corresponding eigenfunctions $e(\unu;\xi)$ with $e(0;\xi)=u_*^\prime$ near $\unu=0$ such that
\beq\label{e:eigen}
\widehat{\caL_*}(\unu)e(\unu;\xi)=\lambda(\unu)e(\unu;\xi),
\eeq
where
\beq
\label{e:eg-exp}
\bcs
&\lambda(\unu)=-d_{\|}\nu_1^2-d_{\perp}|\unu_h|^2+\caO(|\unu|^4),\\
&e(\unu;\xi)=u_*^\prime+\rmi\nu_1u_{*,k}-\nu_1^2e_{21}-|\unu_h|^2e_{2h}+\caO(|\unu|^3),
\ecs
\eeq
with 
\begin{subequations}
\label{e:eg-cof}
\begin{align}
d_{\|}&=\frac{2}{\int_{\T_{2\pi}}(u_*^\prime)^2\rmd \xi}\left[\frac{\rmd }{\rmd k}\left(\int_{\T_{2\pi}}\big((u_*^\dprime)^2-(u_*^\prime)^2\big)\rmd \xi\right)\Bigg|_{k=1}+\int_{\T_{2\pi}}\big(3(u_*^\dprime)^2-(u_*^\prime)^2\big) \rmd \xi\right],\label{e:eg-cof-1}\\
d_{\perp}&=\frac{2}{\int_{\T_{2\pi}}(u_*^\prime)^2\rmd \xi}\left[ \int_{\T_{2\pi}}\big((u_*^\dprime)^2-(u_*^\prime)^2\big)\rmd \xi\right], \label{e:eg-cof-2}\\
e_{21}&=\widehat{\caL_*}(0)^{-1}\bigg( d_{\|}u_*^\prime+4(u_{*,k}^\tprime+u_{*,k}^\prime)+2(3u_*^\tprime+u_*^\prime)\bigg),\ \text{and} \label{e:eg-cof-3}\\
e_{2h}&=\widehat{\caL_*}(0)^{-1}\bigg(d_{\perp}u_*^\prime+2(u_*^\tprime+u_*^\prime)\bigg).\label{e:eg-cof-4}
\end{align}
\end{subequations}
\end{itemize}
\el
\bpf
Without loss of generality, we look for eigenfunctions in the form 
\beq\label{e:decom}
e=u_*^\prime+e_\perp, \text{ with } \langle e_\perp, u_*^\prime \rangle_{L^2(\T_{2\pi})} =0.
\eeq
 Substituting the ansatz \eqref{e:decom} into \eqref{e:eigen}, we obtain
 \beq\label{e:egexp}
\mathcal{F}(e_\perp, \lambda, \unu):= \widehat{\caL_*}(0)e_\perp-\lambda u_*^\prime-\Big[(\widehat{\caL_*}(0)-\widehat{\caL_*}(\unu))(u_*^\prime+e_\perp)+\lambda e_\perp\Big]=0.
 \eeq
Given that $\mathcal{F}(0,0,0)=0$ and that $\partial_{e_\perp, \lambda}\mathcal{F}(0,0,0)$ is an invertible bounded linear operator, an implicit-function-theorem argument shows that,
given $|\unu|$ small, we can solve $\mathcal{F}=0$ for $e_\perp$ and $\lambda$ in terms of $\unu$. 

In order to obtain coefficients in the expansion of $\lambda$ and $e_\perp$,  we first note that $u_*$ is even. We find $\widehat{\caL_*}(\unu)\Big(e(-\nu_1, \unu_h;-\xi)\Big)=\lambda(-\nu_1, \unu_h)e(-\nu_1,\unu_h;-\xi)$, and therefore
\[
\lambda(-\nu_1, \unu_h)=\lambda(\nu_1,\unu_h).
\]
As a consequence, $\lambda=\caO(|\unu|^2)$ only admits even terms in $\unu$.  Next, substituting 
\begin{align*}
\lambda(\unu)=-\sum_{j=1}^nd_{2j}\nu_j^2+\caO(|\unu|^4),\qquad 
e_\perp(\unu;\xi)=\rmi\nu_1e_1-\sum_{j=1}^ne_{2j}\nu_j^2+\caO(|\unu|^3),
\end{align*}
into \eqref{e:egexp}, we find at subsequent orders in $\nu$,
\begin{subequations}
\label{e:cof}
\begin{align}
\rmi\nu_1:\quad \widehat{\caL_*}(0)e_{1}&=4\partial_x(\partial_x^2+1)u_*^\prime=\widehat{\caL_*}(0)u_{*,k}, \label{e:cof1}\\
\nu_1^2: \quad  \widehat{\caL_*}(0)e_{21}&=d_{21}u_*^\prime+4(u_{*,k}^\tprime+u_{*,k}^\prime)+2(3u_*^\tprime+u_*^\prime),\label{e:cof2}\\
\nu_j^2: \quad  \widehat{\caL_*}(0)e_{2j}&=d_{2j}u_*^\prime+2(u_*^\tprime+u_*^\prime),\label{e:cof3}
\end{align}
\end{subequations}
where $2\leq j\leq n$. The last equality in \eqref{e:cof1} is obtained by differentiating \eqref{e:k} with respect to $k$. As a result, $e_1=u_{*,k}$. Solvability of \eqref{e:cof2} implies that the right hand side of the equation should lie in the range of $\widehat{\caL_*}(0)$, that is,
\[
\langle d_{21}u_*^\prime+4(u_{*,k}^\tprime+u_{*,k}^\prime)+2(3u_*^\tprime+u_*^\prime), u_*^\prime \rangle_{L^2(\T_{2\pi})}=0,
\]
implying that $d_{21}=d_{\|}$ and $e_{21}$ takes the form shown in \eqref{e:eg-cof-3}.
A similar analysis on \eqref{e:cof3} shows that for all $2\leq j\leq n$, we have $d_{2j}\equiv d_\perp$ and $e_{2j}=e_{2h}$, which concludes the proof.
\epf
\begin{remark}[Establishing Hypothesis \ref{h:1}]
For $\mu\gtrsim 0$, one can readily show existence of $u_*(\xi;k)$ and establish that the eigenvalue $\lambda=0$ of $\widehat{\caL_*}(0)$ is algebraically simple as assumed in Lemma \ref{l:1.2}, using either Lyapunov-Schmidt or center-manifold reduction. 
Spectral perturbation theory and a careful expansion of the eigenvalue problem for $\widehat{\caL_*}(\unu)$ near $\unu\sim 0$ then also show that $\lambda=0$ does not belong to the spectrum of $\widehat{\caL_*}(\unu)$ for
$\unu\neq 0$ and that $d_{\perp},d_{\|}>0$, as long as 
$k\in K_\mathrm{stab}(\mu)=(k_\mathrm{zz}(\mu),k_\mathrm{eck}(\mu))$, where $k_\mathrm{zz}$ and $k_\mathrm{eck}$ denote the zigzag and Eckhaus stability boundaries; see for instance \cite{mielke_1997}. 
Moreover, $1\in K_\mathrm{stab}(\mu)$ for $\mu\gtrsim 0$, which establishes Hypothesis \ref{h:1} for $\mu\gtrsim 0$.  
\end{remark}
The leading-order expansion $\lambda=-d_\|\nu_1^2-d_\perp |\nu_h|^2$ is the same as the expansion of an anisotropic Laplacian 
\begin{equation}\label{e:deltaeff}
    \Delta_\mathrm{eff}=d_\| \partial_{x_1}^2+d_\perp\Delta_{x_\perp},
\end{equation}
which we refer to as the effective diffusivity near the striped pattern. For $d_\|,d_\perp> 0$, we also define the action of the diffusivity in Fourier space and the scaled norms
\[
D\unu^2:=d_{\|}\nu_1^2+d_\perp |\unu_h|^2, \quad |\ux|_\mathrm{eff}^2:=x_1^2/d_{\|}+|\ux_\perp|^2/d_\perp.
\]
From Remark \ref{r:owendis}(i), we infer that $\Delta_\mathrm{eff}$ is Fredholm on $L^2_\gamma$ for suitable $\gamma$ and infer Fredholm properties. Our goal is next to establish this correspondence for $\caL_*$. 

We start by introducing a pseudo-differential operator with the symbol $\lambda(\nu)$ close to the origin. Recalling the cut-off function $\chi_r$ from \eqref{e:rechi}, we define
\begin{equation}\label{e:Lam}
  \widehat{\LambdaO}(\unu):=\textcolor{blue}{-}D\unu^2+\chi_r(\unu) Q_3(\unu), \qquad Q_3(\nu):=\lambda(\unu)+D\unu^2=\caO(|\unu|^3),
\end{equation}
so that $\widehat{\LambdaO}(\unu)=\lambda(\unu)$ for $|\unu|<r$ and $\widehat{\LambdaO}(\unu)=-D\unu^2$ for $|\unu|>2r$, thus mimicking a pure Laplacian for large and intermediate $\unu$ and the exact behavior of $\caL_*$  near the origin. Since the symbol $\chi_rQ_3$ is compactly supported and smooth, it induces a bounded operator on $L^2$ so that $\LambdaO$, defined through the symbol $\widehat{\LambdaO}$ is a bounded perturbation of the Laplacian on $L^2$. The next lemma states that $\LambdaO-\Delta_\mathrm{eff}$ is small when considered as an operator on Kondratiev spaces.
\begin{lemma}\label{l:konsmall}
 For any $\gamma>0$, there is $r$ sufficiently small such that $\LambdaO$ is a small bounded perturbation of $\Delta_\mathrm{eff}:M^{2,2}_{\gamma-2}\to L^2_\gamma$, that is, for any $\varepsilon>0$,there exists $r_0(\veps)>0$ such that 
 \[
 \|\Delta_\mathrm{eff}-\LambdaO\|\leq \varepsilon \text{ for } r <r_0(\varepsilon).
 \]
 In particular  $\LambdaO$ is Fredholm for $r$ sufficiently small and  $\gamma+n/2\not\in\Z\backslash(2,n)$, with trivial kernel for $\gamma+n/2>2$ and the same index as $\Delta_\mathrm{eff}$.
\end{lemma}
\bpf
It is sufficient to establish the smallness estimate for the Fourier symbol for integer values of $\gamma=k\in\Z^+$, that is, to show that, 
\begin{equation}\label{e:q3}
    \|\chi_r(\unu)Q_3(\unu) u(\unu)\|_{H^{k}}\leq C r
    \left(\|u(\unu)\|_{H^{k-2}}+ 
      \|\unu u(\unu)\|_{H^{k-1}}+
      \sum_{i,j}\|\nu_i\nu_j u(\unu)\|_{H^{k}}
    \right),
\end{equation}
then concluding smallness for fractional $\gamma$ using interpolation.
For this, it is sufficient to consider $Q$ homogeneous of degree 3, with additional factors easily absorbed into $\chi$ in the following argument. Differentiating $\chi_r Q_3 u$ with respect to $\unu$, and grouping terms, one easily finds, for an appropriate constant $C(k)$, 
\begin{align*}
    \|\chi_r(\unu)Q_3(\unu) u(\unu)\|_{H^{k}}\leq &C \left(\|\sum_j \chi_r(\unu)\nu_j \sum_{\ell,m} \left(\nu_\ell\nu_m u(\unu)\right)^{(k)}\|_{L^2}
    + \|\sum_j \chi_r(\unu)\nu_j\sum_{\ell} \left(\nu_\ell u(\unu)\right)^{(k-1)}\|_{L^2} \right. \\
   &\left. + \|\sum_j \chi_r(\unu)\nu_j u(\unu)^{(k-2)}\|_{L^2}    + \sum_{3\leq j\leq k}\|\chi_r(\unu)u(\unu)^{(k-j)}\|_{L^2}\right). 
\end{align*}
The first three terms are readily estimated using that $\chi_r(\unu)\unu$ is small in $L^\infty$ for $r$ small, while the last term can be estimated using the fact that $\|u^{(k-3)}\|_{L^p(\Omega)}\leq \|u^{(k-3)}\|_{H^1(\Omega)}$ for some $p>2$ and compact domain $\Omega$ and then estimating $\|\chi_r(\unu)u^{(k-3)}\|_{L^2}$ using H\"older's inequality, which proves \eqref{e:q3} and the lemma. 
\epf

\subsection{Fredholm properties of $\caL_*$ via conjugacies
}\label{s:3.2}

We introduce the shorthand  notation for operators 
\[
\caL_0:=\caL_*(1-\caL_*)^{-1}, \quad 
\caL_\Lambda:=\LambdaO(1-\LambdaO)^{-1},
\]
both closed and densely defined on $L^2_\gamma$,
with  Bloch-Fourier decompositions
\[
\caB\circ\caL_0\circ\caB^{-1}=\int_{\Omega}\widehat{\caL_0}(\unu)\rmd \unu, \quad \caB\circ\caL_\Lambda\circ\caB^{-1}=\int_{\Omega}\widehat{\caL_\Lambda}(\unu)\rmd \unu.
\]
We start our analysis by examining these decompositions on algebraically weighted spaces $L^2_\gamma(\R^n)$ for any $\gamma\geq0$. For $\gamma>0$, $L^2_\gamma(\R^n)$ is a subspace of $L^2(\R^n)$ and thus $\caB$ defined on $L^2_\gamma(\R^n)$ is simply the restriction of $\caB$ to \eqref{e:BFT} on $L^2_\gamma(\R^n)$. Our first goal is to characterize the range or $\caB$ restricted to $L^2_\gamma$. Fourier transform is an isomorphism between $L^2_\gamma(\R^n)$ and $H^\gamma(\R^n)$, defined for instance using interpolation between $[\gamma]:=\max\{k\in\Z,k\leq \gamma\}$ and $[\gamma]+1$. Introducing
\beq\label{e:gamma-}
\gamma_-:=\{k\in\Z\mid  \gamma-n/2-1\leq k<\gamma -n/2\},
\eeq
as the largest integer that is \textit{strictly} less than $\gamma-n/2$,
we have the continuous embedding
\[
 H^\gamma(\R^n) \hookrightarrow C^{\gamma_-, \tau}(\R^n) , \qquad \text{ for } \gamma_-\in\Z,\ \tau\in (0,1],\ \gamma_-+\tau+n/2\leq\gamma,
\]
\bl\label{l:blochimage}
The range  of $\caB$ restricted to $L^2_\gamma$, $\gamma\geq 0$, contains functions which extend to twist-periodic functions that are of class $H^\gamma$, that is, 
\[
\caB(L^2_\gamma(\R^{n}))=\left\{ u_B(\nu_1,\unu_\perp;\xi)\mid u_B(\nu_1,\unu_\perp;\xi)=\rme^{-\rmi \nu_1\xi} U(\nu_1,\unu_\perp;\xi), U\in H^\gamma_\mathrm{per}(\Omega, L^2(\T_{2\pi}))\right\},
\] 
and $U(\nu_1,\unu_\perp;\xi)=U(\nu_1+1,\unu_\perp;\xi)$. 
We write 
$
\caB (L^2_\gamma)=: H^\gamma_\mathrm{tw}
.$
\el
\bpf 
The lemma is a direct consequence of the definition of the Bloch wave transform. 
We have, writing $\widehat{u}$ for the Fourier transform,
\[
u_B(\nu_1,\unu_\perp;\xi)=\sum_\ell e^{i\ell_1 \xi}\widehat{u}(\ell_1+\nu_1,\unu_\perp),\qquad 
U(\nu_1,\unu_\perp;\xi)=\sum_\ell e^{i(\ell_1+\nu_1) \xi}\widehat{u}(\ell_1+\nu_1,\unu_\perp)
.\]
Choosing now $\gamma=0$, we find that norms of $u_B$ and $U$ in $L^2(\Omega,L^2(\T_{2\pi}))$ are equivalent to the norm of $\widehat{u}$ or $u$ in $L^2(\R^n)$.  Also, $U$ is  periodic in $\nu_1$ which can be seen by shifting the index of summation. If $u\in L^2_1(\R^n)$, we have $\widehat{u}\in H^1(\R^n)$ 
so that $\nabla_{\unu}u_B,\nabla_{\unu_\perp} U \in L^2(\Omega,L^2(\T_{2\pi}))$. 
Differentiating the equation for $U$ with respect to $\nu_1$, we also find that $\partial_{\nu_1}U\in  L^2(\Omega,L^2(\T_{2\pi}))$ after noticing that the multiplicative term $\rmi\xi U$ is bounded. One then readily finds that the norm of $U\in H^1(\Omega,L^2(\T_{2\pi}))$ is equivalent to the norm of $u$ in $L^2_1$. Higher derivatives are obtained in the same fashion and fractional  values of $\gamma$ by interpolation.
\epf
From the lemma, we have equivalent norms on  $ H^\gamma_\mathrm{tw}$, either induced by $\caB$, or by the isomorphism with $H^\gamma_\mathrm{per}(\Omega, L^2(\T_{2\pi})$,
\begin{equation}\label{e:gammaeq}
\|v\|_{H^\gamma_\mathrm{tw}}:=\|\caB^{-1}v\|_{L^2_\gamma(\R^n)}\sim \|\rme^{\rmi\nu_1\xi }v\|_{H^\gamma_\mathrm{per}(\Omega,L^2(\T_{2\pi}))},
\end{equation}
where the latter norm is defined via Fourier series.

To analyze $\caL_0$ and $\caL_\Lambda$ in terms of their Bloch-Fourier counterparts, we  construct a $C^\infty$-isomorphism $\iota:=\int_{\Omega}\iota(\unu)\rmd\unu$ with
\beq
\label{e:iota}
\iota(\unu)=\iota_0(\unu)\oplus\iota_\perp(\unu): L^2(\T_{2\pi})\longrightarrow L^2(\T_{2\pi}).
\eeq
Here, the maps $\iota_0,\iota_\perp$ are constructed as isomorphisms, 
\[
\iota_0: \mathrm{span}\{1\} \to\mathrm{span}\{\widetilde{e}(\unu)\}, \quad \iota_\perp: \mathrm{span}\{ 1 \}^\perp\to \mathrm{span}\{\widetilde{e}(\unu)\}^\perp,
\]
where $\widetilde{e}(\unu):= e_0+\chi_{4r_0}(\unu)[e(\unu)-e_0]$ and where $\chi_{4r_0}$ is a smoothed characteristic function for a ball of radius $4r_0$, as defined in \eqref{e:chi}. Therefore set
$\iota_0(\unu)(1)= \widetilde{e}(\unu)$, and  define $\iota_\perp(\unu)$ as the identity on $\mathrm{span}\{ 1,\widetilde{e}(\unu)\}^\perp$ and mapping $\mathrm{span}\{ 1\}^\perp\cap \mathrm{span}\{ 1,\widetilde{e}(\unu)\}$  to $\mathrm{span}\{ \widetilde{e}(\unu)\}^\perp\cap \mathrm{span}\{ 1,\widetilde{e}(\unu)\}$ preserving norm.
%
%
%
We later will restrict to modes in balls of radius $2r_0$ and $r_0$, hence the peculiar choice of $4r_0$. 


To summarize, there exists an operator $\iota$ such that $\iota$ and $\iota^{-1}$ are $C^\infty$-smooth as a mapping from $\Omega$ to $\mathscr{B}(L^2(T_{2\pi}))$. Moreover, $\iota, \iota^{-1}\in W^{k,\infty}\left(\Omega, \mathscr{B}(L^2(T_{2\pi}))\right)$, for any $k\in\Z^+$, and
\[
\widecheck{\iota}:= \caB\circ \iota\circ \caB^{-1}: L^2_\gamma(\R^n)\to L^2_\gamma(\R^n)
\]
is an isomorphism for any $\gamma\in\N$, which, can be generalized to $\gamma\geq0$ via interpolation. Note that $\iota$ trivially preserves the twist periodicity since it is simply constant with respect to $\unu$ near the boundary $\partial\Omega$. We state this observation in the following lemma.
\bl\label{l:iota}
For any $\gamma\geq0$, the operator $\widecheck{\iota}: L^2_\gamma(\R^n)\to L^2_\gamma(\R^n)$ is an isomorphism.
\el
We next introduce a Bloch-Fourier multiplier $M:=\int_{\Omega}M(\unu)\rmd\unu$ through
\beq
\label{e:BF-M}
M(\unu)\circ\widehat{\caL_\Lambda}(\unu)=\iota^{-1}(\unu) \circ\widehat{\caL_0}(\unu)\circ\iota(\unu),
\eeq
with counterpart in  physical space, $\caM:= \caB^{-1}\circ M\circ\caB$. 

Our goal is to show that $\caM$, restricted to the closure of $ \rg(\caL_\Lambda)$ is an isomorphism, which will be the key step to establishing that $\caL_\Lambda$ and $\caL_0$ share the same Fredholm properties. First notice that 
\begin{equation}\label{e:conj}
\caM\circ \caL_\Lambda \circ \widecheck{\iota^{-1}}=\widecheck{\iota^{-1}}\circ\caL_0,
\end{equation}
implies that it is sufficient  to show that $\caM$ restricted to the closure of $ \rg(\caL_\Lambda)$ is an isomorphism.

To start with, we summarize some relevant properties of $M$ and its inverse $M^{-1}$ in the following lemma.
\bl\label{l:02}
The operator $M(\unu): L^2(\T_{2\pi})\to L^2(\T_{2\pi})$ and its inverse are of class  $C^{\infty}\left(\Omega, \mathscr{B}(L^2(T_{2\pi}))\right)$, and uniformly bounded on $\Omega$.
\el
\bpf 
The results are straightforward when $\unu$ is away from zero, using \eqref{e:BF-M} and the fact that both $
\widehat{\caL_\Lambda}(\unu)$ and $\widehat{\caL_0}(\unu)$ are bounded invertible. For $|\unu|$ sufficiently small (say, $|\unu|<r_0$), the operators $\widehat{\caL_\Lambda}$ and $\iota^{-1}\circ\widehat{\caL_0}\circ\iota$ on $\langle 1\rangle\oplus \langle 1\rangle^\perp$ admit the decomposition
\[
\widehat{\caL_\Lambda}\cong \begin{pmatrix}  \lambda(\unu) (1-\lambda(\unu))^{-1}& 0\\ 0 & \widehat{\caL_\Lambda}\mid_{\langle 1\rangle^\perp}\end{pmatrix}, \qquad
\iota^{-1}\circ\widehat{\caL_0}\circ\iota\cong 
\begin{pmatrix} \lambda(\unu) (1-\lambda(\unu))^{-1}& 0\\ 0 & \iota_\perp^{-1}\circ\widehat{\caL_0}\circ\iota_\perp
\end{pmatrix}, 
\]
and thus the operator $M(\unu)$ is given, in this decomposition, through
\[
M(\unu)\cong \begin{pmatrix} 1 & 0\\ 0 & M_{11}(\unu)\end{pmatrix},
\]
where  $M_{11}= \iota_\perp^{-1}\circ\widehat{\caL_0}\circ\iota_\perp \circ  \left(\widehat{\caL_\Lambda}\mid_{\langle 1\rangle^\perp}\right)^{-1}$. It is readily seen that $M(\unu)$ is  smooth and  uniformly bounded invertible for $|\unu|<r_0$.
\epf
We note that conjugation with the effective Laplacian $\Delta_\mathrm{eff}
$ would give regularity $W^{2,\infty}$ of the symbol $M(\unu)$, thus limiting the ranges of allowed $\gamma$-values. Smoothness of $M$, Lemma \ref{l:02}, and of $\iota$, Lemma \ref{l:iota}, imply that \eqref{e:conj} defines a conjugacy of the operators $\mathcal{L}_*$ and $\LambdaO$, thus guaranteeing Fredholm properties of $\caL_*$. 

We now readily obtain the following result on  Fredholm properties for the linearization at stripes analogous to Theorem \ref{t:McOwen}, roughly showing that up to the cokernel of pseudo-harmonic polynomials, the operator is bounded invertible.
\bt\label{thm:fredholm}
For space dimension $n\geq 2$ and spatial weight $n/2<\gamma+n/2\not\in \Z\backslash(2,n)$, the linearization at stripes, considered as a closed operator, 
\beq\label{e:lin}
\caL_*: \mathcal{D}_\gamma(\caL_*)\subset H^4_{\gamma-2}(\R^n)\to L^2_\gamma(\R^n)
\eeq
is Fredholm. We have the characterization of kernel and cokernel as follows: 
\begin{enumerate}
\item for $\gamma+n/2\in(2, n)$, \eqref{e:lin} is an isomorphism.
\item for $\gamma+n/2\in(n+m, n+m+1)$ and $m\in\N$,  \eqref{e:lin} is an injection with its closed range equal to 
$
R_\gamma= \caH_{\leq m}^{\perp_\gamma}(\caL_*).
$
\end{enumerate}
\et
\bpf
The result is an immediate consequence of the smoothness of $\iota$ and $M$, which provide a smooth conjugacy of Fourier symbols. 
\epf


Counting dimensions of the kernel and exploiting the conjugacy with the Laplacian, one can then completely characterize the kernel of $\caL_*$ when restricted to polynomially growing functions. We collect here only the relevant information for at most linear growth, and refer to the appendix for a more comprehensive study. 

\bl\label{l:cokernel}
The kernel of the linear operator $\caL_*$ in spaces of functions with at most linear growth is spanned by the  {\rm pseudo-harmonic polynomials}, $H_0,H_{1,j},$ $j=1\,\ldots,n$. 
\el

\bpf By smallness of the perturbation in Kondratiev spaces and conjugacy, the kernel is at most $n+1$-dimensional. A direct calculation readily demonstrates that the $n+1$ functions $H_0,H_{1,j}$ are linearly independent and belong to the kernel. 
\epf

We may then construct a projection  $\mathcal{Q}$ onto the cokernel of $\caL_*$, using the fact that the range is $L^2$-orthogonal to the kernel in $L^2_{-\gamma}$. We will only use this projection in the specific case where cokernel has dimension $n+1$ and write  
\begin{equation}\label{e:caQ}\mathcal{Q}R_\gamma=0,\qquad  
\mathrm{dim}\,(\mathrm{Rg}\,\mathcal{Q})=\mathrm{dim}\,(\mathrm{coker}\,\caL_*).
\end{equation}
Note that this projection is not uniquely defined as there is no canonical basis for a complement of the range in $L^2_\gamma$. 

\subsection{Mode filters}\label{s:3.3}
Unfortunately, the loss of localization predicted in Theorem \ref{thm:fredholm} does not allow us to close a nonlinear argument when quadratic nonlinearities are present, as is the case in our situation.  We therefore develop here and in the next section more refined mapping properties of  $\caL_*$. Roughly speaking, the operator should be bounded invertible for Bloch modes away from the neutral mode, without loss of localization. For close-to-neutral modes, we expect to lose 2 degrees of localization as predicted in Theorem \ref{thm:fredholm}. Localization is regained  by applying operators obtained by conjugating derivatives with $\iota$, since applying derivatives after inverting $\Lambda_\mathrm{eff}$ gains derivatives. 
Splitting off neutral modes will be accomplished by the mode filters that we introduce below. We then introduce spaces that characterize gain of localization via derivatives. 

We define the mode filters first in Fourier and then in physical space, through
\beq\label{e:Pr}
\caP_r\caB f(\unu;\xi):= \chi_{r}(\unu)\frac{\langle \caB f(\unu;\cdot), e(\unu;\cdot)\rangle_{L^2(\T_{2\pi})}}{\| e(\unu;\cdot)\|_{L^2(\T_{2\pi})}^2}e(\unu; \xi), \qquad P_{0,r}:=\caB^{-1}\caP_r\caB,\quad  P_{h,r}:=\id-P_{0,r},
\eeq
where $\chi_{r}$ is the smoothed characteristic function of a small neighborhood of the origin in $\R^n$ defined in \eqref{e:chi}, and $0<r\leq 2r_0$, where $r_0$ was used in the definition of the conjugacy $\iota$ in \eqref{e:iota}. We will throughout use the fact  that
\begin{equation}\label{e:mfs}
    P_{0,r}P_{h,2r}=P_{h,2r}P_{0,r}=0, \text{ for any }r\in(0,r_0],
\end{equation}
which is immediate from  $\chi_r(1-\chi_{2r})\equiv0$.
Note that $P_{0,r}$ and $P_{h,r}$ would be projections if $\chi_{r}$ was an actual characteristic function. Nevertheless, $P_{0,r}$ and $P_{h,r}$ commute with the operator $\caL_*$. More precisely, we have, for any $u\in\caD_\gamma(\caL_*)$, 
\beq
\label{e:LP=PL}
\widehat{\caL_*}\caP_r\caB u=\lambda(\unu)\chi_r(\unu)\frac{\langle \caB u(\unu;\cdot), e(\unu;\cdot)\rangle_{L^2(\T_{2\pi})}}{\| e(\unu;\cdot)\|_{L^2(\T_{2\pi})}^2}e(\unu; \xi)=\caP_r\widehat{\caL_*}\caB u.
\eeq

A natural consequence of the commutativity property is that $P_{0,r}$ and $P_{h,r}$ are well-defined as mappings on the range $R_\gamma$; that is,
$
P_{0,r} (R_\gamma)\subseteq R_\gamma$ and $P_{h,r} (R_\gamma)\subseteq R_\gamma$.
Note that $P_{0,r}$ and $P_{h,r}$ are defined and bounded  on $L^2_\gamma$ for any $\gamma$, based on the smoothness of $\chi_{r}$. We can therefore decompose any function in the range, 
\begin{equation}\label{e:fdec}
f=\caL_*u\in R_\gamma=\caH^{\perp_\gamma}_{\leq m}(\caL_*) \implies f=\overbrace{P_{0,r}f}^{=:f_0}+\overbrace{P_{h,r}f}^{=:f_h}, \quad f_0,f_h\in R_\gamma.
\end{equation}
From now on, we assume that the kernel of $\caL_*$ is trivial, $\gamma+n/2>2$. Statements below will hold true for other values of $\gamma$ when taking norms after quotienting the kernel. 

Boundedness of $P_{0/h,r}$ and Theorem  \ref{thm:fredholm} give that, for $\max\{2,n/2\}<\gamma+n/2\not\in\Z\backslash(2,n)$,
\begin{equation}\label{e:1stest}
\|u_0\|_{H^4_{\gamma-2}}\leq C \|f_0\|_\gamma, \quad \|u_h\|_{H^4_{\gamma-2}}\leq C \|f_h\|_\gamma,\qquad \text{where }\  u_0:=\caL_*^{-1}f_0, \  u_h:=\caL_*^{-1}f_h, 
\end{equation}
and we note that $u_0= P_{0,r}u$ since $\caL_*$ and $P_{0,r}$ commute.
Our next goal is to refine those estimates with a further decomposition of the neutral component $u_0$ to characterize the gain of localization and a Bloch-Fourier characterization of the stable component $u_h$.

A more explicit expression of mode filters can be gained by expanding the mode filter operator $\caP_r$ in $\unu$ in Bloch-Fourier space. In particular, keeping only the leading-order term $u_*^\prime$ of the eigenfunction $e(\unu;\xi)$, we define the simplified mode filter,
\beq
\label{e:tildePr}
\widetilde{\caP}_r\caB u:=\chi_{r}(\unu)\frac{\langle (\caB u)(\unu;\cdot), e(\unu;\cdot)\rangle_{L^2(\T_{2\pi})}}{\| e(\unu;\cdot)\|_{L^2(\T_{2\pi})}^2}u_*^\prime(\xi).
\eeq
Transforming back to physical space, the multiplication by $u_*'$ factors, so that
\beq
\caB^{-1}\widetilde{\caP}_r\caB u=\caB^{-1}\left[\widetilde{\caP}_r\caB u\right]=\caB^{-1}\left[\chi_{r}(\unu)\frac{\langle (\caB u)(\unu;\cdot), e(\unu;\cdot)\rangle_{L^2(\T_{2\pi})}}{\| e(\unu;\cdot)\|_{L^2(\T_{2\pi})}^2}\right] u_*'(x_1)=:
\left[\widetilde{P}_{0,r} u\right]u_*'(x_1).
\eeq
Note that the last expression is a pointwise product in $\ux$ of the envelope $\widetilde{P}_{0,r} u$ and the periodic function $u_*'(x_1)$.
We then define the complementary simplified  mode filter as 
\beq
\widetilde{P}_{h,r}u:=(\id-u_*^\prime(x_1)\widetilde{P}_{0,r})u.
\eeq
Summarizing, we have 
\beq\label{e:2nddec}
u=\overbrace{P_{0,r}u}^{u_0}+\overbrace{P_{h,r}u}^{u_h}=u_*^\prime(x_1)\overbrace{\widetilde{P}_{0,r}u}^{=:\widetilde{u}_0}+\overbrace{\widetilde{P}_{h,r}u}^{=:\widetilde{u}_h},
\eeq
where we introduce the notation $\widetilde{u}_0$ for the envelope and $\widetilde{u}_h$ for the remainder.

In order to state properties of our decomposition, we also introduce the  modified Kondratiev spaces, 
\begin{equation}\label{e:dk}
  \|u\|_{M^{m,k,p}_\gamma(\R^n)}=\left(\sum_{0\leq|\underline{\alpha}|\leq m} \|D^{\underline{\alpha}} u  \|^p_{L^p_{\min(\gamma+k,\gamma+|\underline{\alpha}|)}}\right)^{1/p},
\end{equation}
where  $|\underline{\alpha}|=\alpha_1+\ldots+\alpha_n$,  $D^{\underline{\alpha}} $ denotes partial derivatives with multi-index convention. The following result gives the desired refined version of  \eqref{e:1stest}.
\begin{proposition}
\label{p:kondrat}
For $u\in \caD_\gamma(\caL_*)$ with  $n\geq 2$, $\max\{2,n/2\}<\gamma+n/2\not\in\Z\backslash(2,n)$, with the decomposition $\caL_*u=f=f_0+f_h\in R_\gamma$ as in \eqref{e:fdec},  and the decomposition $u=u_*^\prime\widetilde{u}_0+\widetilde{u}_h$ as in \eqref{e:2nddec}, there exists a constant $C>0$ so that 
\[
\|\widetilde{u}_0\|_{M^{4,2,2}_{\gamma-2}}\leq C \|f\|_{L^2_\gamma},\qquad
\|\widetilde{u}_h\|_{H^4_{\gamma-1}}\leq C  \|f\|_{L^2_\gamma}.
\]
In other words, we can recover the loss of localization by taking derivatives in the envelope, $\widetilde{u}_0$, and lose just one degree of localization in the complement, $\widetilde{u}_h$. 
\end{proposition}
We prepare the proof with the following characterization of the neutral component. 
\begin{lemma}\label{l:rgf}
For $u\in \caD_\gamma(\caL_*)$ with $n\geq 2$, $\max\{2,n/2\}<\gamma+n/2\not\in\Z\backslash(2,n)$,  there exists $w\in M^{2,2}_{\gamma-2}(\R^n)$ such that 
\beq
\label{e:tildeu0}
\caB \widetilde{u}_0=g(\unu)\widehat{w}(\unu),
\quad  
\text{ with  }\quad g(\unu)=\frac{1}{2\pi}\chi_{r}(\unu).
 \eeq
\end{lemma}
\begin{Proof}
For $u\in \caD_\gamma(\caL_*)$, we define $f=\caL_*u$, and consider the following commutative diagram,
\begin{equation}
    \begin{tikzpicture}
  \matrix (m) [matrix of math nodes,row sep=6em,column sep=6em,minimum width=4em]
  {
     L^2_{\gamma-2}(\R^n) & H^2_{\gamma-2}(\R^n)\hookleftarrow M^{2,2}_{\gamma-2}(\R^n) &L^2_{\gamma}(\R^n) & L^2_{\gamma}(\R^n)\\
L^2_{\gamma-2}(\R^n)& H^4_{\gamma-2}(\R^n)\hookleftarrow \caD_\gamma(\caL_*) & & L^2_{\gamma}(\R^n)\\};
  \path[-stealth]
    (m-1-1) edge node [left] {$\widecheck{\iota}$} (m-2-1)
        edge node [above] {$(1-\Lambda_{{\rm eff}})^{-1}$} (m-1-2)
        (m-1-2) edge node [above] {$\Lambda_{{\rm eff}}$} (m-1-3)
    (m-1-3) edge node [above] {$\caM$} (m-1-4)
     (m-1-4) edge node [right] {$\widecheck{\iota}$} (m-2-4)
     (m-2-1) edge node [above] {$(1-\caL_*)^{-1}$} (m-2-2)
(m-2-2) edge node [above] {$\caL_*$} (m-2-4);
\end{tikzpicture}
\label{t:cd}
\end{equation}
From the above diagram, together with Theorem \ref{t:McOwen} and Theorem \ref{thm:fredholm}, there exist a unique $w\in M^{2,2}_{\gamma-2}(\R^n)$, such that
\[
\caB u=\widehat{\caL_*}^{-1}\circ\iota\circ M\circ \widehat{\Lambda_{{\rm eff}}} (\caB w).
\]
This yields
\[
\begin{aligned}
\caB \widetilde{u}_0&=
\chi_{r}(\unu)\frac{\langle (\caB u)(\unu;\cdot), e(\unu;\cdot)\rangle_{L^2(\T_{2\pi})}}{\| e(\unu;\cdot)\|_{L^2(\T_{2\pi})}^2}\\
&=\frac{\chi_{r}(\unu)}{\| e(\unu;\cdot)\|_{L^2(\T_{2\pi})}^2}\langle \widehat{\caL_*}^{-1}\circ\iota(\unu)\circ M(\unu)\circ \widehat{\Lambda_{{\rm eff}}}(\unu)( \caB w) (\unu;\cdot), e(\unu;\cdot)\rangle_{L^2(\T_{2\pi})}\\
&=\frac{\chi_{r}(\unu)}{2\pi\lambda(\unu)}\langle M(\unu)\circ \widehat{\Lambda_{{\rm eff}}}(\unu)( \caB w) (\unu;\cdot), 1\rangle_{L^2(\T_{2\pi})}\\
&=\frac{\chi_{r}(\unu)}{2\pi\lambda(\unu)}\langle \lambda(\unu)( \caB w) (\unu;\cdot), 1\rangle_{L^2(\T_{2\pi})}\\
&=\frac{\chi_{r}(\unu)}{2\pi} \langle (\caB w)(\unu;\cdot), 1 \rangle_{L^2(\T_{2\pi})}\\
&=\frac{\chi_{r}(\unu)}{2\pi} \widehat{w}(\unu),\\
\end{aligned} 
\]
which concludes the proof.
\end{Proof}

\begin{Proof}[ of Proposition \ref{p:kondrat}.] We first prove the estimate on $\widetilde{u}_0$. We recall from \eqref{e:tildeu0} in Lemma \ref{l:rgf} that
$
\caB \widetilde{u}_0=g(\unu)\widehat{w}(\unu).
$
We need to verify that applying up to 2 derivatives improves localization by one degree for each derivative,  while the action of higher derivatives
is simply bounded. First, according to the commutative diagram \eqref{t:cd}, we have
\beq\label{e:h0-1}
w=\Lambda_\mathrm{eff}^{-1}\circ\caM^{-1}\circ\widecheck{\iota}^{-1}f, \quad \text{with } \|w\|_{M^{2,2}_{\gamma-2}}\leq C\|f\|_{L^2_\gamma}
\eeq
where the estimate is derived from the mapping properties of $\Lambda_\mathrm{eff}$, Theorem \ref{t:McOwen} and Lemma \ref{l:konsmall}. Next, for $|\underline{\alpha}|\leq 2$, we have 
\beq\label{e:h0-2}
\begin{aligned}
\|D^{\underline{\alpha}}\widetilde{u}_0\|_{L^2_{\gamma-2+|\underline{\alpha}|}}
= &\|\caB (D^{\underline{\alpha}}
\widetilde{u}_0)
\|_{H^{\gamma-2+|\underline{\alpha}|}(\Omega,L^2(\T_{2\pi}))}\\
= &\|(\partial_\xi+\rmi \nu_1)^{\alpha_1}(\rmi\unu_\perp)^{\underline{\alpha}_\perp}\caB \widetilde{u}_0\|_{H^{\gamma-2+|\underline{\alpha}|}(\Omega,L^2(\T_{2\pi}))}\\
= &\|(\partial_\xi+\rmi \nu_1)^{\alpha_1}(\rmi\unu_\perp)^{\underline{\alpha}_\perp}g\widehat{w}\|_{H^{\gamma-2+|\underline{\alpha}|}(\Omega,L^2(\T_{2\pi}))}\\
= &\|(\rmi \nu_1)^{\alpha_1}(\rmi\unu_\perp)^{\underline{\alpha}_\perp}g\widehat{w}\|_{H^{\gamma-2+|\underline{\alpha}|}(\Omega,L^2(\T_{2\pi}))}\\
\leq &C\|(\rmi \unu)^{\underline{\alpha}}\widehat{w}\|_{H^{\gamma-2+|\underline{\alpha}|}(\R^n)}\\
= &C\|D^{\underline{\alpha}}w\|_{L^2_{\gamma-2+|\underline{\alpha}|}(\R^n)}.
\end{aligned}
\eeq
Noting that any higher derivative of $\widetilde{u}_0$ can be applied to $w$ as shown in the above inequality \eqref{e:h0-2} and the rest can be absorbed into the constant, we conclude that, for any $|\underline{\alpha}|>2$, 
\beq\label{e:h0-3}
\begin{aligned}
\|D^{\underline{\alpha}}\widetilde{u}_0\|_{L^2_{\gamma}}
\leq &C\|w\|_{M^{2,2}_{\gamma-2}(\R^n)}. 
\end{aligned}
\eeq
Combining \eqref{e:h0-1}, \eqref{e:h0-2} and \eqref{e:h0-3}, we have
\[
\|\widetilde{u}_0\|_{M^{4,2,2}_{\gamma-2}}\leq C\|f\|_{L^2_\gamma}.
\]
We are left to show the estimate for $\widetilde{u}_h$. We recall from \eqref{e:2nddec} that
\beq\label{e:tildeuh}
\widetilde{u}_h=u_h+\overbrace{u_0-u_*^\prime(x_1)\widetilde{u}_0}^{=:Q_ru_0},
\eeq
and note that 
\[
(\caB Q_ru_0)(\unu;\xi)=g(\unu)\widehat{w}(\unu)(e(\unu;\xi)-u_*^\prime(\xi))=g(\unu)\widehat{w}(\unu)(\rmi\nu_1 e_1(\unu;\xi)+\caO(|\unu|^2)),
\]
where we apply the same techniques as in the inequality \eqref{e:h0-2} and readily conclude that 
\beq
\label{e:rhoest}
\|Q_ru_0\|_{H^4_{\gamma-1}}\leq C\|w\|_{M^{2,2}_{\gamma-2}}\leq C\|f\|_{L^2_\gamma}.
\eeq
Due to the fact that $f_h$ does not admit any close-to-neutral modes, the action of $\caL_*^{-1}$ on $f_h$ in Bloch-space is bounded invertible so that we readily obtain that 
\beq
\label{e:uhest}
\|u_h\|_{H^4_\gamma(\R^n)}=\|\caL_*^{-1}f_h\|_{H^4_\gamma(\R^n)}\leq C\|f_h\|_{L^2_\gamma(\R^n)}.
\eeq
Combing \eqref{e:tildeuh}-\eqref{e:uhest} concludes the proof.
\end{Proof}

\subsection{Refined mapping properties of $\caL_*$}\label{s:3.4}
We show refined estimates based on the mode filters and modulation decomposition introduced above. Recall that
\begin{equation}
\label{e:tilde-P-0h}
    \widetilde{P}_{0, r}u=\caB^{-1}\left(\frac{ \langle\caB(P_{0,r}u)(\unu,\cdot), e (\unu,\cdot) \rangle_{L^2(\T_{2\pi})}}{\|e(\unu;\cdot)\|_{L^2(\T_{2\pi})}^2} \right), \quad \widetilde{P}_{h,r}u=u-u_*^\prime \widetilde{P}_{0, r}u.
\end{equation}
Roughly speaking, we show that the linear operator loses 2 degrees of localization \emph{within} the neutral component, no localization in a complement, and has cross terms between neutral component and complement that lose 1 degree of localization. 
\bp\label{p:II} 
Let $r\in(0,r_0)$, $n\geq 2$, and $\max\{2,n/2\}<\gamma+n/2\not\in\Z\backslash(2,n)$. Given $f\in L^2_{\gamma-1}$ such that $P_{0,2r}f\in R_\gamma$, we have $f\in R_{\gamma-1}$ and
\begin{equation}\label{e:lin_est}
\|\widetilde{P}_{0, r} \caL_*^{-1}f\|_{M^{4,2,2}_{\gamma-2}}+
\|\widetilde{P}_{h,r} \caL_*^{-1}f\|_{H^4_{\gamma-1}}\leq C\left(\|f\|_{L^2_{\gamma-1}}+\|P_{0,2r}f\|_{L^2_{\gamma}}\right).
\end{equation} 
\ep

 \bpf
We claim that if $P_{0,2r}f\in R_\gamma$, then $f\in R_{\gamma-1}$. Thanks to Theorem \ref{thm:fredholm}, the claim is obviously true for $\gamma<1+n/2$; for $\gamma>1+n/2$, we have that for any $H\in \caH_{\leq (\gamma-1)-}(\caL_*)$,
\[
\langle H, f \rangle_{L^2(\R)}=\langle H, P_{0,2r}f \rangle_{L^2(\R)}=0,
\]
since $\langle H, P_{0,2r}f \rangle_{L^2(\R)}$ is a linear combination of coefficients in the Taylor expansion of 
\[
\langle \caP_{2r}\caB  f, e(\unu;\cdot)\rangle_{L^2(\T_{2\pi})}=\chi_{2r}(\unu)\langle \caB  f, e(\unu;\cdot)\rangle_{L^2(\T_{2\pi})},
\]
at $\unu=0$.
Recalling the equality in \eqref{e:mfs} and the fact that $P_{h,2r}$ and $\caL_*$ are commutative, we conclude that
\[
P_{0,r}\caL_*^{-1}P_{h,2r}f
 =P_{0,r}P_{h,2r}\caL_*^{-1}f =0.
 \]
From this, together with the definition of $\widetilde{P}_{0,r}$ in \eqref{e:tilde-P-0h} and the fact that  $Q_r=P_{0,r}-\widetilde{P}_{0,r}$ as in \eqref{e:2nddec} and \eqref{e:tildeuh}, we readily find that
\[
\widetilde{P}_{0,r}\caL_*^{-1}P_{h,2r}f=0, \qquad Q_r\caL_*^{-1}P_{h,2r}f=(P_{0,r}-u^\prime(x_1)\widetilde{P}_{0,r})\caL_*^{-1}P_{h,2r}f=0,
  \]
 which, in combination with Proposition \ref{p:kondrat}, yields 
 \beq\label{e:P0Linv}
 \begin{aligned}
\|\widetilde{P}_{0,r}\caL_*^{-1}f\|_{M^{4,2,2}_{\gamma-2}}&=\|\widetilde{P}_{0,r}\caL_*^{-1}P_{0,2r}f\|_{M^{4,2,2}_{\gamma-2}}\leq C\|P_{0,2r}f\|_{L^2_\gamma},\\
  \|Q_r\caL_*^{-1}f\|_{H^4_{\gamma-1}}&=\|Q_r\caL_*^{-1}P_{0,2r}f\|_{H^4_{\gamma-1}}\leq C\|P_{0,2r}f\|_{L^2_\gamma}.
 \end{aligned}
 \eeq
 On the other hand, we have 
 \[
P_{h,r}\caL_*^{-1}f=\caL_*^{-1}P_{h,r}f,
 \]
 where, again, $P_{h,r}f$ admits only stable modes of $\caL_*$ and thus we readily conclude that
 \beq\label{e:PhLinv}
 \|P_{h,r}\caL_*^{-1}f\|_{H^4_{\gamma-1}}=\|\caL_*^{-1}P_{h,r}f\|_{H^4_{\gamma-1}}\leq C\|f\|_{L^2_{\gamma-1}}.
 \eeq
 Combining \eqref{e:P0Linv}--\eqref{e:PhLinv})and noting that $\widetilde{P}_{h,r}=Q_r+P_{h,r}$, we find the desired estimate \eqref{e:lin_est}.
 \epf

\section{Modulation ansatz and fixed point scheme}\label{s:4}

We substitute the  ansatz  
\beq\label{e:ansatz}
u=u_\psi + w=u_\psi+u_*^\prime \widetilde{w}_0+\widetilde{w}_h,
\eeq
into the stationary Swift-Hohenberg equation, where 
\[
u_\psi =u_*(\psi(\ux); |\nabla \psi(\ux)|), 
\]
with phase modulation  $\psi(\ux)=x_1+\Theta(\ux)$ and $\Theta(\ux)$ small, localized, to be determined. Recall that we focus on the physically relevant cases $n=2,3$. We incorporate the uniform phase shift $a_0$ into $g$, writing
\beq \label{e:shan}
-(\Delta+1)^2(u_\psi+w)+\mu (u_\psi +w) - (u_\psi+w)^3+\veps g(x_1-a_0, \ux_\perp)=0.
\eeq
Now define 
\beq
\label{e:IFT-nt}
\begin{matrix}
\caL_*:=-(\Delta+1)^2+\mu-3u_*^2, &\quad\qquad\qquad\caN(\psi, w):=3(u_\psi^2-u_*^2)w+3u_\psi w^2+w^3, \\
\caR(\psi):=-[-(\Delta+1)^2+\mu]u_\psi +u_\psi^3,   & \caG(\vphi,\veps):=-\veps g(x_1-a_0,\ux_\perp),
\end{matrix}
\eeq
formally apply $\caL_*^{-1}$ to \eqref{e:shan}, and expand to find
\beq\label{e:ift-form}
w=\caL_*^{-1}\left(\caR(\psi) + \caN(\psi, w)+\caG(a_0, \veps)\right),
\eeq
which, equivalently, can be considered as an equation in $(\widetilde{w}_0,\widetilde{w}_h)\in X_{\gamma-2}:=M^{4,2,2}_{\gamma-2}\times H^4_{\gamma-1}$,
\beq\label{e:ift-form-ref}
\begin{pmatrix}
\widetilde{w}_0 \\ \widetilde{w}_h
\end{pmatrix}=\begin{pmatrix}
\widetilde{P}_{0,r} \\
\widetilde{P}_{h,r}
 \end{pmatrix} \caL_*^{-1}\left(\caR(\psi) + \caN(\psi, \widetilde{w}_0 u_*^\prime +\widetilde{w}_h)+\caG(a_0, \veps)\right)
.\eeq
Clearly, \eqref{e:ift-form} possesses a trivial solution $w=0$ at $\veps=0$. The strategy then is, in a somewhat straightforward fashion, to establish contraction mapping properties for the right-hand side when $\veps$ is small. In fact, once one establishes that the right-hand side is well defined, one finds differentiability and can solve with the Implicit Function Theorem. There are however two key difficulties: the linear operator $\caL_*$ is, first, not onto and, second,  loses localization when inverting. The first difficulty is addressed by the introduction of the phase modulation $\psi$, which adds degrees of freedom that compensate for the cokernel. The second difficulty is addressed through the refined estimates \eqref{e:lin_est}. Without these refined estimates or the phase modulation $\psi$, it appears to be impossible to close such a nonlinear iteration argument. 

It turns out that the nonlinear argument needs to be carried out in spaces with localization that enforces a cokernel spanned by linearly growing pseudo-harmonic polynomials. This can be compensated by smoothed versions of the Green's function and its (discrete) derivatives. The Green's function itself does however introduce weak decay in the far field, or even logarithmic growth in two dimensions, leading to solutions that are not observed. The bounded element of the cokernel is therefore better compensated by an appropriate translate $a_0$. Since those translates enter the equation only at order $\veps$, we introduce a dummy term $\alpha \phi(\ux)$ at first, to compensate for the constants in the cokernel. Given a solution in terms of $\veps$, we then solve $\alpha=0$ with variable $a_0$ after obtaining expansions for $\alpha(\veps)$ and dividing by $\veps=0$.

Slightly more specifically, we change \eqref{e:ift-form-ref} to
\beq\label{e:ift-form-ref2}
0=-\begin{pmatrix}
\widetilde{w}_0 \\ \widetilde{w}_h
\end{pmatrix}+\begin{pmatrix}
\widetilde{P}_{0,r} \\
\widetilde{P}_{h,r}
 \end{pmatrix} \caL_*^{-1}\left(\caR(\psi(\underline{a})) + \caN(\psi(\underline{a}), \widetilde{w}_0 u_*^\prime +\widetilde{w}_h)+\caG(a_0, \veps)+\alpha \phi\right),
\eeq
and solve as a nonlinear equation
\begin{equation}\label{e:fnl}
    0=\mathcal{F}(\widetilde{w}_0,\widetilde{w}_h,\underline{a},a_0,\alpha;\veps).
\end{equation}
Here, $(\widetilde{w}_0,\widetilde{w}_h)\in X_{\gamma-2}$, $a_0,\alpha\in\R$, and $\underline{a}\in\R^n$, where $n$ is the dimension of ambient space and $\psi(\underline{a})$ is linear in $\underline{a}$, consisting of a parameterization of phase corrections by the cokernel defined in the next section. Using Proposition \ref{p:II}, well-posedness of \eqref{e:fnl} requires 
\beq\label{e:nec}
\caN+\caG+\caR+\alpha\phi\in R_{\gamma-1}, \qquad P_{0, 2r}(\caN+\caG+\caR+\alpha\phi)\in R_{\gamma},
\eeq
when $(\widetilde{w}_0,\widetilde{w}_h)\in X_{\gamma-2}$. Note that the inclusions \eqref{e:nec} contain statements on \emph{localization} and,  \emph{cokernel conditions}. Satisfying both conditions and showing that the resulting iteration scheme defines a contraction mapping is then accomplished in the following steps:
\begin{enumerate}
    \item 
        construct $\psi(\ua)$, linear in $\ua$, and $\phi$ so that $\caR(\psi(\cdot))$ and $\alpha \phi$ span the cokernel of $\caL_*$ and $\caR\in L^2_{\gamma-1}$, with $\ua\in \R^n$ and $\alpha\in\R$;
    \item 
        show that $\caN$ is continuously differentiable  on $X_{\gamma-2}\times\R^n $ into $L^2_{\gamma-1}$ with vanishing derivative at the origin;
    \item 
        show that $P_{0, 2r}\caN$ is well defined and continuously differentiable from $X_{\gamma-2}\times\R^n $ into $L^2_\gamma$; 
    \item 
        show that $\caR+\caN+\caG+\alpha\phi\in R_{\gamma-1}$, $P_{0, 2r}(\caN+\caG+\caR+\alpha\phi)\in R_{\gamma}$, that is, both are  in the range, for a choice $\ua=\Psi_a(\widetilde{w}_0,\widetilde{w}_h,a_0)$, $\alpha=\Psi_\alpha(\widetilde{w}_0,\widetilde{w}_h,a_0)$ with $\Psi_a:X_{\gamma-2}\to\R^n$, $\Psi_{a_0}:X_{\gamma-2}\to\R$ smooth with vanishing derivative at the origin;
    \item 
       substitute the expressions for $\ua$ and $\alpha$ into  \eqref{e:ift-form-ref2}, solve for $(\widetilde{w}_0,\widetilde{w}_h)$ with the implicit function theorem as functions of $a_0$, and  expand the solution with $\alpha=\alpha(a_0;\veps)$;
    \item expand $\alpha=\veps\alpha_1({a_0})+\rmO(\veps^2)$ and find $\alpha=0$ at a nondegenerate zero of $\alpha_1$ using Hypothesis~\ref{h:2}.
\end{enumerate}
These steps will be carries out in the subsequent sections. 

\subsection{Modulations by multi-pole modulations}

We recall first that
\[
\caR(\psi)=[-(\Delta+1)^2+\mu]u_\psi-u_\psi^3=-(\Delta^2+2\Delta)u_\psi+ |\nabla\psi|^4u_\psi^\qprime+2|\nabla\psi|^2u_\psi^\dprime,
\]
where $u_\psi=u(\psi; |\nabla\psi|)$, and we used the equation for $u$ as a periodic solution.  We introduce the ansatz
\begin{equation}\label{e:psi}
\Theta(\ux):=(1-\chi_1(\ux))\ua\cdot\nabla G(\ux),
\end{equation}
 and exploit the phase modulation
\[
\psi(\ux)=x_1+\Theta(\ux).
\]
In order to evaluate decay, we first derive expressions for derivatives, 
{\allowdisplaybreaks
\begin{align}\label{e:lapex}
\nabla u =&u^\prime \nabla \psi + u_{k}\nabla|\nabla\psi|, \nonumber\\
\Delta u=&\nabla\cdot(u^\prime\nabla\psi)+\nabla\cdot(u_k\nabla|\nabla\psi|) \nonumber\\
=&u^\dprime|\nabla\psi|^2+2u^\prime_k\nabla\psi\cdot\nabla|\nabla\psi|+u_{kk}\left|\nabla|\nabla\psi|\right|^2+u^\prime\Delta\psi+u_{k}\Delta|\nabla\psi|, \nonumber\\
\Delta^2u=&\left[(\Delta u^\dprime)|\nabla\psi|^2+2(\nabla u^\dprime)\cdot(\nabla |\nabla\psi|^2)+u^\dprime(\Delta|\nabla\psi|^2)\right]+ \nonumber\\
&2\left\{(\Delta u^\prime_k)\nabla\psi\cdot\nabla|\nabla\psi|+2\left(\nabla u^\prime_k\right)\cdot \left[\nabla\left(\nabla\psi\cdot\nabla|\nabla\psi|\right)\right]+u_{k}^\prime\Delta\left(\nabla\psi\cdot\nabla|\nabla\psi|\right)\right\}+  \nonumber\\
&\left[(\Delta u_{kk})\left|\nabla|\nabla\psi|\right|^2+2(\nabla u_{kk})\cdot(\nabla \left|\nabla|\nabla\psi|\right|^2)+u_{kk}(\Delta\left|\nabla|\nabla\psi|\right|^2)\right]+\nonumber\\
&\left[(\Delta u^\prime) (\Delta\psi)+2(\nabla u^\prime)\cdot\nabla\left(\Delta\psi\right)+u^\prime\Delta^2\psi\right]+\\
&\left[(\Delta u_{k}) (\Delta|\nabla\psi|)+2(\nabla u_{k})\cdot\nabla\left(\Delta|\nabla\psi|\right)+u_{k}\Delta^2|\nabla\psi|\right]\nonumber\\
=&u^\qprime|\nabla\psi|^4+4u_{k}^\tprime |\nabla\psi|^2\left(\nabla\psi\cdot \nabla|\nabla\psi|\right)+2u^\dprime_{kk}\left[|\nabla\psi|^2\left|\nabla|\nabla\psi|\right|^2+2\left(\nabla\psi\cdot \nabla|\nabla\psi|\right)^2\right]+\nonumber\\
&4u_{kkk}^\prime\left|\nabla|\nabla\psi|\right|^2\left(\nabla\psi\cdot \nabla|\nabla\psi|\right)+u_{kkkk}\left|\nabla|\nabla\psi|\right|^4+ 2u^\tprime\nabla\cdot\left(|\nabla\psi|^2 \nabla\psi\right)+\nonumber\\
&2u_{k}^\dprime\nabla\cdot\left[|\nabla\psi|^2\nabla|\nabla\psi|+2\left(\nabla\psi\right)\cdot\left(\nabla|\nabla\psi|\right)\nabla\psi\right]+\nonumber\\
&2u^\prime_{kk}\nabla\cdot\left[\left|\nabla|\nabla\psi|\right|^2\nabla\psi+2\left(\nabla\psi\right)\cdot\left(\nabla|\nabla\psi|\right)\nabla|\nabla\psi|\right]+ 2u_{kkk}\nabla\cdot\left(\left|\nabla|\nabla\psi|\right|^2 \nabla|\nabla\psi|\right)+ \nonumber\\
&u^\dprime[\frac{1}{2}\Delta^2(\psi^2)-\psi\Delta^2\psi]+u^\prime_k[\Delta^2(\psi|\nabla\psi|)-|\nabla\psi|\Delta^2\psi-\psi\Delta^2|\nabla\psi|]+\nonumber\\
&u_{kk}[\frac{1}{2}\Delta^2(|\nabla\psi|^2)-|\nabla\psi|\Delta^2|\nabla\psi|]+u^\prime\Delta^2\psi+u_{k}\Delta^2|\nabla\psi|. \nonumber
\end{align}
}
Since the gradient of the Green's function is homogeneous of degree $1-n$, we express decay in orders of algebraic decay, that is, we write $\caO(m)$ when a function is bounded by $|\ux|^{-m}$ up to a constant for $|\ux|$ sufficiently large. 

For large $\ux$, we have
$
\psi(\ux)=x_1+\caO(n-1),
$
and 
\beq\label{e:est}
\psi_{x_1}=1+\caO(n), \quad \nabla_{\ux_\perp}\psi=\caO(n), \quad \nabla|\nabla\psi|=\nabla\psi_{x_1}+\caO(n+1). 
\eeq
We now expand $\caR$ as
\beq\label{e:carexp}
\begin{aligned}
\caR(\psi)=&-4u^\tprime_k |\nabla\psi|^2\left(\nabla\psi\cdot \nabla|\nabla\psi|\right)-2u^\tprime\nabla\cdot\left(|\nabla\psi|^2 \nabla\psi\right)-\\
&2u^\dprime_k\nabla\cdot\left[|\nabla\psi|^2\nabla|\nabla\psi|+2\left(\nabla\psi\right)\cdot\left(\nabla|\nabla\psi|\right)\nabla\psi\right]-u^\dprime[\frac{1}{2}\Delta^2(\psi^2)-\psi\Delta^2\psi]-\\
&u^\prime_k[\Delta^2(\psi|\nabla\psi|)-|\nabla\psi|\Delta^2\psi-\psi\Delta^2|\nabla\psi|]-u^\prime\Delta^2\psi-u_{k}\Delta^2|\nabla\psi|-\\
&2[2u^\prime_k\nabla\psi\cdot\nabla|\nabla\psi|+u^\prime\Delta\psi+u_{k}\Delta|\nabla\psi|]+\caO(n+3)\\
=&-4u_{k}^\tprime\Delta_{x_1}\psi-2u^\tprime(3\Delta_{x_1}\psi+\Delta_{\ux_\perp}\psi)-2u^\dprime_k(3\Delta_{x_1}\psi+\Delta_{\ux_\perp}\psi)_{x_1}-4u^\dprime(\Delta\psi)_{x_1}-\\
&4u^\prime_k\Delta_{x_1}(\Delta\psi)-u^\prime\Delta^2\psi-u_{k}(\Delta^2\psi)_{x_1}-\\
&4u^\prime_k\Delta_{x_1}\psi-2(u^\prime\Delta\psi+u_{k}(\Delta\psi)_{x_1} )+\caO(n+3)\\
=&(-4u^\tprime_k-6u^\tprime-4u^\prime_k-2u^\prime)\Delta_{x_1}\psi+(-2u_*^\tprime-2u^\prime)\Delta_{\ux_\perp}\psi+\caO(n+2).
\end{aligned} 
\eeq
With this calculation, we are ready to conclude that $\caR$ contributes terms to the fixed point iteration with sufficient localization.
\bl \label{l:rloc}
For $\gamma+n/2<2+n$ we have 
\beq\label{e:Rinc}
\caR\in L^2_{\gamma-1},\qquad P_{0,2r} \caR\in L^2_\gamma.
\eeq
Moreover, $\caR$ is smooth as a function of $\ua$.
\el
\bpf
Since $\caR=\caO(n+1)$, we find that $\caR\in L^2_{\gamma-1}$, which gives the first inclusion in \eqref{e:Rinc}. In order to find the extra localization after applying the mode filter; that is, to see that $P_{0, 2r}\caR\in L^2_\gamma(\R^n)$, we write
\[
D_{\|}:=-4u^\tprime_k-6u^\tprime-4u^\prime_k-2u^\prime, \qquad D_{\perp}:=-2u_*^\tprime-2u^\prime,
\]
and exploit the expansion of $\caR$,
\[
P_{0,2r}\caR=P_{0,2r}\left( D_{\|}\psi_{x_1x_1}+D_{\perp}\Delta_{\ux_\perp}\psi\right)+\caO(n+2).
\]
Noting that any term of the order $\caO(n+2)\in L^2_\gamma(\R^n)$ for any $\gamma<2+n/2$ and that the function $d_{\|}\psi_{x_1x_1}+d_\perp \Delta_{\ux_\perp}\psi$ has compact support, we are left to show that 
\[
P_{0,2r}\bigg( \left(D_{\|}\psi_{x_1x_1}+D_\perp\Delta_{\ux_\perp}\psi\right)-\left( d_{\|}\psi_{x_1x_1}+d_\perp \Delta_{\ux_\perp}\psi\right)u_*^\prime\bigg)\in L^2_\gamma(\R^n),
\] 
or, equivalently,
\[
\caB P_{0,2r}\bigg(\left(D_{\|}\psi_{x_1x_1}+D_\perp\Delta_{\ux_\perp}\psi\right)-\left(d_{\|}\psi_{x_1x_1}+d_\perp\Delta_{\ux_\perp}\psi\right)u_*^\prime\bigg)\in H^\gamma_\mathrm{tw}.
\]
We therefore expand
\[
\begin{aligned}
&\caB P_{0,2r}\bigg(\left(D_{\|}\psi_{x_1x_1}+D_\perp\Delta_{\ux_\perp}\psi\right)-\left(d_{\|}\psi_{x_1x_1}+d_\perp\Delta_{\ux_\perp}\psi\right)u_*^\prime\bigg)\\
=&\left\langle \left(D_{\|}\caB\psi_{x_1x_1}+D_\perp\caB\Delta_{\ux_\perp}\psi\right)-\left(d_{\|}\caB\psi_{x_1x_1}+d_\perp\caB\Delta_{\ux_\perp}\psi\right)u_*^\prime, e\right\rangle_{L^2(T_{2\pi})}\frac{\chi_{2r}(\unu)}{\|e(\unu;\cdot)\|^2_{L^2(\T_{2\pi})}}e(\unu;\xi)\\
=&\left\langle \left(D_{\|}\caB\psi_{x_1x_1}+D_\perp\caB\Delta_{\ux_\perp}\psi\right)-\left(d_{\|}\caB\psi_{x_1x_1}+d_\perp\caB\Delta_{\ux_\perp}\psi\right)u_*^\prime, u_*^\prime+\caO(|\unu|)\right\rangle_{L^2(T_{2\pi})}\frac{\chi_{2r}(\unu)}{\|e(\unu;\cdot)\|^2_{L^2(\T_{2\pi})}}e(\unu;\xi).
\end{aligned}
\]
We again could apply a discrete derivative argument on the higher order terms and show that they are at least of order $\caO(n+2)$ and thus in $L^2_\gamma$.  We now note that
\[
d_{\|}=\frac{\langle D_{\|},u_*^\prime \rangle_{L^2(\T(2\pi))}}{\|u_*^\prime\|^2_{L^2(\T_{2\pi})}}, \qquad d_\perp=\frac{\langle D_\perp,u_*^\prime \rangle_{L^2(\T(2\pi))}}{\|u_*^\prime\|^2_{L^2(\T_{2\pi})}},
\]
and rewrite the essential part of the leading order term,
\[
\begin{aligned}
&\left\langle \left(D_{\|}\caB\psi_{x_1x_1}+D_2\caB\Delta_{\ux_\perp}\psi\right)-\left(d_{\|}\caB\psi_{x_1x_1}+d_\perp\caB\Delta_{\ux_\perp}\psi\right)u_*^\prime, u_*^\prime\right\rangle_{L^2(T_{2\pi})}\\
&=\left\langle \caB\psi_{x_1x_1},(D_{\|}-d_{\|}u_*^\prime)u_*^\prime\right\rangle_{L^2(T_{2\pi})}+\left\langle \caB\Delta_{\ux_\perp}\psi,(D_\perp-d_\perp u_*^\prime)u_*^\prime\right\rangle_{L^2(T_{2\pi})}
\end{aligned}
\]
where the two functions $(D_{\|}-d_{\|}u_*^\prime)u^\prime$ and $(D_\perp-d_\perp u_*^\prime)u^\prime$ are $2\pi$-periodic with zero average. We may therefore integrate  by parts and find that the derivative with respect to $\xi$ of $\caB\psi_{x_1x_1}$ and $\caB\Delta_{\ux_\perp}\psi$ yields an additional degree of decay in physical space. 

Smoothness of $\caR$ and $P_{0,2r}\caR$ as functions of $\ua$ is readily obtained as follows. Dependence on $\ua$ through $\psi$ is linear so that boundedness implies smoothness. Dependence on $\ua$ through $u$ yields additional degrees of localization upon differentiation, and smooth dependence readily follows. 
%
\epf
Finally, we also introduce for later use  the correction 
\begin{equation}\label{e:phi}
    \phi(\ux)=\chi_R(\ux)u_*'(x_1),
\end{equation}
with an arbitrary, fixed $R>0$.

\subsection{Localization and smoothness of $\caN$}

In order to evaluate superposition operators in weighted spaces, we prepare with embedding results that recover sharp localization from weighted Sobolev and Kondratiev spaces.
\bl\label{l:l-infty}
There is a constant  $C$ such that for any $u_h\in H^4_\gamma$, we have
\[
|u_h(\ux)|\leq C\|u_h\|_{H^4_{\gamma}(\R^n)}\langle\ux\rangle^{-\gamma}.
\]
\el
\bpf
This is immediate using the embedding of $H^4$ into $L^\infty$ for the function $u(\ux)\langle\ux\rangle^\gamma$.
\epf

\bl
\label{l:alg-dom}
There is a continuous embedding of $M^{2,2,2}_{\gamma+n/p-n/2}(\R^n)$ into $L^p_{\gamma}(\R^n)$, for $n=2,3$, $1\leq p\leq\infty$; that is,
\[
\begin{cases}
    \|u_0\|_{L^p_{\gamma}}\leq C\|u_0\|_{M^{2,2,2}_{\gamma+n/p-n/2}}, & \quad 1\leq p <\infty; \\ 
    \|u_0\|_{L^\infty_{\gamma}}\leq C\|u_0\|_{M^{2,2,2}_{\gamma-n/2}}, & \quad p=\infty.
\end{cases}
\]
\el
\bpf 
We prove the global weighted embedding via local classical Sobolev embedding on an exponential annulus decomposition of $\R^n$. More specifically, we introduce a family of annuli,
\[
A_j=\begin{cases}
\{\ux\mid 2^j\leq |\ux|< 2^{j+1}\}, \quad j\in\N;\\
\{\ux\mid |\ux|<1\}, \quad j=-1.
\end{cases}
\]
On each annulus $A_j$ ($j\geq 0$), we introduce the change of variable,
\[
u_j(\uy)=u(2^{j}\uy),\quad \uy\in A_0,
\]
so that, under the change of coordinates, the function $u$ restricted to each $A_j$ ($j\geq 0$) is now defined on the same annulus $A_0$, yielding with universal changing constants $C$,
\[
\begin{aligned}
\|u\|_{L^p_{\gamma}}^p
=&\sum_{j=-1}^\infty\int_{A_j}\langle\ux\rangle^{p\gamma}|u(\ux)|^p\rmd \ux
=\int_{A_{-1}}\langle\ux\rangle^{p\gamma}|u(\ux)|^p\rmd \ux+ \sum_{j=0}^\infty\int_{A_j}\langle\ux\rangle^{p\gamma}|u(\ux)|^p\rmd \ux\\
\leq &C\left[\|u\|_{H^2(A_{-1})}^p+\sum_{j=0}^\infty 2^{[p\gamma+n]j}\int_{A_0}|u_j(\uy)|^p\rmd \uy\right]
=
C\left[\|u\|_{H^2(A_{-1})}^p+\sum_{j=0}^\infty 2^{[p\gamma+n]j}\|u_j\|_{L^p(A_0)}^p\right]\\
\leq &C\left[\|u\|_{H^2(A_{-1})}^p+\sum_{j=0}^\infty 2^{[p\gamma+n]j}\|u_j\|_{H^2(A_0)}^p\right],
\end{aligned}
\]
by Sobolev embedding in dimension $n\leq 3$. Next, exploiting the Kondratiev scaling, 
{\small
\[
\begin{aligned}
\|u\|_{L^p_{\gamma}}^p
\leq &C\left[\|u\|_{H^2(A_{-1})}^p+\sum_{j=0}^\infty 2^{[p\gamma+n]j}\left(\int_{A_0}|u_j(\uy)|^2\rmd\uy+ \int_{A_0}|\nabla_{\uy}u_j|^2\rmd\uy+ \int_{A_0}|\nabla^2_{\uy}u_j|^2\rmd\uy\right)^{p/2}\right]\\
\leq &C\left[\|u\|_{H^2(A_{-1})}^p+\sum_{j=0}^\infty 2^{[p\gamma+n]j}\left(2^{-nj}\int_{A_j}|u(\ux)|^2\rmd\ux+2^{(2-n)j}\int_{A_j}|\nabla_{\ux}u(\ux)|^2\rmd\ux+\right.\right.\\
& \left.\left.+2^{(4-n)j}\int_{A_j}|\nabla^2_{\ux}u(\ux)|^2\rmd\ux\right)^{p/2}\right]\\
\leq &C\left[\|u\|_{H^2(A_{-1})}^p+\left(\sum_{j=0}^\infty 2^{[2\gamma+2n/p-n]j}\int_{A_j}|u(\ux)|^2\rmd\ux\right)^{p/2}+ \left(\sum_{j=0}^\infty 2^{[2\gamma+2n/p-n+2]j}\int_{A_j}|\nabla_{\ux}u(\ux)|^2\rmd\ux\right)^{p/2}\right.\\
& \left.
+\left(\sum_{j=0}^\infty 2^{[2\gamma+2n/p-n+4]j}\int_{A_j}|\nabla^2_{\ux}u(\ux)|^2\rmd\ux\right)^{p/2}\right]\\
\leq & C\|u\|_{M^{2,2,2}_{\gamma+n/p-n/2}}^p.
\end{aligned}
\]}
The case $p=\infty$ can be treated in the same fashion.
\epf

We are now ready to evaluate the superposition operators. 
We first establish mapping properties of  
$
\caN
$.

\bl\label{l:N}
For $n=2,3$, and $\gamma+n/2\in (n+1, n+2)$, we have that $\caN$, defined through 
\begin{align*}
\caN: X_{\gamma-2}\times\R^n&\to L^2_{\gamma-1},\\
(\widetilde{w}_0,\widetilde{w}_h,\ua) &\mapsto 3(u_*^2-u_\psi^2)w-3u_*w^2+3(u_*-u_\psi)w^2-w^3,
\end{align*}
with 
$w=\widetilde{w}_0+\widetilde{w}_h$, $\psi$ as in \eqref{e:psi} is well-defined and smooth on a small neighborhood of the origin, with  vanishing derivative at the origin.  
\el

\bpf
We decompose
\[
\caN=\overbrace{3(u_*^2-u_\psi^2)w}^{=:\caN_1}-\overbrace{3u_*w^2}^{=:\caN_2}+\overbrace{3(u_*-u_\psi)w^2-w^3}^{=:\caN_3}.
\]
Recalling that 
\[
\psi=x_1+\Theta(\ux)=x_1+(1-\chi(\ux))\ua\cdot\nabla G(\ux),
\]
we note that  $\caN_1$ is bounded by $\Theta w$. Since $\Theta$ is $\caO(1-n)$ and $w\in L^2_{\gamma-2}$, this shows $\caN_1\in L^2_{\gamma-1}$, depending smoothly on the parameters $\ua$. 
%
To show smoothness and boundedness of of $\caN_2$ and $\caN_3$, we use that $w\in L^2_{\gamma-2}$ and, by Lemma \ref{l:alg-dom}, $w\in L^\infty_{\gamma-2+n/2}$, so that again $w^2\in L^2_{\gamma-1}$, smoothly depending on $w$ as a bilinear bounded map. Derivatives in $\ua$ introduce additional localization and are treated in the same fashion. Since $u_*=u_\psi$ are $\ua=0$, $\caN$ is quadratic at the origin with vanishing derivative.  
\epf


\subsection{Localization and smoothness of $P_{0, 2r}\caN$}
The next lemma establishes smoothness similar to Lemma \ref{l:N}, but with one additional degree of localization.
\bl
For $n=2,3$ and $\gamma+n/2\in (n+1,n+2)$, we have that $P_{0,2r}\caN$, defined through 
\begin{align*}
P_{0,2r}\caN: X_{\gamma-2}\times\R^n&\to L^2_{\gamma},\\
(\widetilde{w}_0,\widetilde{w}_h,\ua) &\mapsto P_{0,2r}\left(3(u_*^2-u_\psi^2)w-3u_*w^2+3(u_*-u_\psi)w^2-w^3\right),
\end{align*}
with 
$w=\widetilde{w}_0+\widetilde{w}_h$, $\psi$ as in \eqref{e:psi}, is well-defined and smooth on a small neighborhood of the origin,  with  vanishing derivative at the origin. 
\el 
\begin{Proof}
We check that $P_{0,2r}\caN=P_{0,2r}(\caN_1+\caN_2+\caN_3)\in L^2_\gamma$ term by term.

\emph{$\mathbf{P_{0,2r}\caN_3\in L^2_\gamma}$.} We just need to show that
\[
\caN_3=3(u_*-u_\psi)w^2-w^3\in L^2_\gamma.
\]
Noting that $w^2\in L^2_{\gamma-1}$ and $u_*-u_\psi=\caO(\Theta)=\caO(1-n)$, we readily conclude that 
$
3(u_*-u_\psi)w^2\in L^2_\gamma.$
To treat $w^3$, we recall from Lemma \ref{l:alg-dom} and Lemma \ref{l:l-infty} that
\[
\widetilde{w}_0(\ux)=\caO\Big(-(\gamma-2+n/2)\Big), \qquad
\widetilde{w}_h(\ux)=\caO\Big(-(\gamma-1)\Big).
\]
The first inequality implies that $\tilde{w}_0^3\in L^2_\gamma$, given that $2(\gamma-2+n/2)\geq2$, that is; $\gamma+n/2\geq 3$, which is true for our choice of $\gamma$ and $n$.
Similarly, we conclude that $\widetilde{w}_0^2\widetilde{w}_h, \widetilde{w}_0\widetilde{w}_h^2, \widetilde{w}_h^3\in L^2_\gamma$ observing that $\widetilde{w}_h$ is stronger localized than $\widetilde{w}_0$. In summary, we established that 
$
w^3\in L^2_\gamma
$ as claimed. 

\paragraph{$\mathbf{P_{0,2r}\caN_2\in L^2_\gamma}$.} We expand $w^2=\widetilde{w}_h^2+2u_*^\prime\widetilde{w}_0\widetilde{w}_h+(u_*^\prime)^2\widetilde{w}_0^2$. Using the inequalities
$
\gamma-1\geq1$ and $\gamma-2+n/2\geq 1$, we readily conclude that
$\widetilde{w}_h^2, \widetilde{w}_0\widetilde{w}_h\in L^2_\gamma$ so that we only need to show that
\[
P_{0,2r}(u_*(u_*^\prime)^2\widetilde{w}_0^2)\stackrel{!}{\in} L^2_\gamma,
\]
or equivalently,
\begin{equation}
\caB P_{0,2r}(u_*(u_*^\prime)^2\widetilde{w}_0^2)=\chi_{2r}(\unu)\frac{\langle u_*(u_*^\prime)^2 \caB \widetilde{w}_0^2(\unu;\cdot), e(\unu;\cdot)\rangle_{L^2(\T_{2\pi})}}{\| e(\unu;\cdot)\|_{L^2(\T_{2\pi})}^2}e(\unu; \xi)\stackrel{!}{\in}  H^\gamma_\mathrm{tw}.\label{e:tbp}
\end{equation}
Inserting the leading order expansion $e(\unu;
\xi)=u_*^\prime(\xi)+\caO(|\unu|)$ into the above expression, we have 
\[
\caB P_{0,2r}(u_*(u_*^\prime)^2\widetilde{w}_0^2)=\langle  \caB \widetilde{w}_0^2(\unu;\cdot), u_*(u_*^\prime)^3+\caO(|\unu|)\rangle_{L^2(\T_{2\pi})}
\frac{\chi_{2r}(\unu)}{\| e(\unu;\cdot)\|_{L^2(\T_{2\pi})}^2}e(\unu; \xi).
\]
In order to take advantage of the $\caO(\unu)$ term which gives additional smoothness in Fourier-Bloch space and extra decay in physical space, we use the discrete derivative 
\[
\delta_j u(\ux):=u(\ux+2\pi\mathbf{e}_j)-u(\ux)=\int_0^{2\pi}\partial_{x_j}u(\ux+t\mathbf{e}_j)\rmd t, \text{ for all }\ux\in\R^n,
\]
where $\mathbf{e}_j$ is the $j$-th unit vector. We find that 
\[
\|\delta_ju\|_{M^{4,1,2}_{\gamma-1}}\leq C\|u\|_{M^{4,2,2}_{\gamma-2}},
\]
and 
\[
\caB \delta_j u= (\rme^{-2\pi\rmi\nu_j}-1)\caB u=\caO(|\nu_j|)\caB u.
\]
As a result, in order to establish \eqref{e:tbp}, 
we need to show that $\delta_j(\widetilde{w}_0^2)\in L^2_\gamma$. Using discrete integration by parts, it then suffices to show that $\widetilde{w}_0\delta_j(\widetilde{w}_0)\in L^2_\gamma$. Noting that $\widetilde{w}_0\in M^{4,2,2}_{\gamma-2}$, we have $\delta_j(\widetilde{w}_0)\in L^2_{\gamma-1}$. As a result, we can conclude from the inequality 
$
\gamma-2+n/2\geq 1,
$
that \[
\langle  \caB \widetilde{w}_0^2(\unu;\cdot), \caO(|\unu|)\rangle_{L^2(\T_{2\pi})}
\frac{\chi_{2r}(\unu)}{\| e(\unu;\cdot)\|_{L^2(\T_{2\pi})}^2}e(\unu; \xi)\in H^\gamma_\mathrm{tw}.
\]
Lastly, given that $u_*(u_*^\prime)^3$ is an odd and $2\pi$-periodic function, it admits a $2\pi$-periodic anti-derivative, denoted as $U(\xi)$. Integration by part then gives 
\[
\begin{aligned}
\langle  \caB \widetilde{w}_0^2(\unu;\cdot), u_*(u_*^\prime)^3\rangle_{L^2(\T_{2\pi})}
&=-\langle  \partial_\xi\caB \widetilde{w}_0^2(\unu;\cdot), U\rangle_{L^2(\T_{2\pi})}\\
&=-\langle  \sum_{k\in\Z}\rmi k\widecheck{\tilde{w}}_0^2(k+\nu_1,\unu_\perp)\rme^{\rmi k\xi}, U\rangle_{L^2(\T_{2\pi})}\\
&=-\langle  \sum_{k\in\Z}\rmi (k+\nu_1-\nu_1)\widecheck{\tilde{w}}_0^2(k+\nu_1,\unu_\perp)\rme^{\rmi k\xi}, U\rangle_{L^2(\T_{2\pi})}\\
&=-\langle  \sum_{k\in\Z}\rmi (k+\nu_1)\widecheck{\tilde{w}}_0^2(k+\nu_1,\unu_\perp)\rme^{\rmi k\xi}, U\rangle_{L^2(\T_{2\pi})}+\langle  \caB\widetilde{w}_0^2(\unu;\cdot), \caO(|\unu|)\rangle_{L^2(\T_{2\pi})}\\
&=-\langle  \caB(2\widetilde{w}_0\partial_{x_1}\widetilde{w}_0)(\unu;\cdot), U\rangle_{L^2(\T_{2\pi})}+\langle  \caB\widetilde{w}_0^2(\unu;\cdot), \caO(|\unu|)\rangle_{L^2(\T_{2\pi})},\\
\end{aligned}
\]
from which  we readily conclude that
\[
\langle  \caB \widetilde{w}_0^2(\unu;\cdot), u_*(u_*^\prime)^3)\rangle_{L^2(\T_{2\pi})}
\frac{\chi_{2r}(\unu)}{\| e(\unu;\cdot)\|_{L^2(\T_{2\pi})}^2}e(\unu; \xi)\in H^\gamma_\mathrm{tw},
\]
as claimed.
\paragraph{$\mathbf{P_{0,2r}\caN_1\in L^2_\gamma}$.} We start with  the decomposition,
\[
\caN_1=3(u_*^2-u_\psi^2)w=-6u_*u_*^\prime \Theta(u_*^\prime\widetilde{w}_0+\widetilde{w}_h)+\caO(|\Theta^2|)(u_*^\prime\widetilde{w}_0+\widetilde{w}_h).
\]
Clearly, we have $\Theta\widetilde{w}_h\in L^2_\gamma$ and $\caO(|\Theta|^2)\widetilde{w}_0\in L^2_\gamma$, thanks to the fact that $\caO(\Theta)=\caO(1-n)$, $\widetilde{\omega}_0\in L^2_{\gamma-2}$, and $\widetilde{w}_h\in L^2_{\gamma-1}$.
We are left to show that 
\[
P_{0,2r}(u_*(u_*^\prime)^2\Theta \widetilde{w}_0)\in L^2_\gamma,
\]
the  proof of which is very similar to the above proof for $P_{0,2r}(u_*(u_*^\prime)^2\widetilde{w}_0^2)$ and omitted here.

\paragraph{Smoothness.} Differentiability properties follow readily from the above estimate after noticing that the map can be viewed as a symmetric, bounded,  multilinear map with bounds given by the previously established estimates. Derivatives with respect to $\psi$ induce additional localization and the derivative at the origin vanishes. 
\end{Proof}
\begin{remark} The inhomogeneity is well defined and smooth by our assumptions on localization. Moreover, we have $P_{0,2r}\caG\in L^2_\gamma(\R^n)$ as a straightforward consequence of $P_{0,2r}$ being bounded and linear on $L^2_\gamma(\R^n)$. 
\end{remark}
\subsection{Projecting onto the range through modulation}

We wish to choose $\ua$ and $\alpha$ so that $P_{0,2r}(\caN+\caG+\caR+\alpha\phi)\in \caR_\gamma$, that is, 
\begin{equation}\label{e:rg}
h_0=\int_{\R^n} H_0 P_{0,2r}(\caN+\caG+\caR+\alpha\phi)\rmd \ux\stackrel{!}{=}0, \qquad
h_\ell=\int_{\R^n} H_{1,\ell}P_{0,2r}(\caN+\caG+\caR+\alpha\phi)\rmd \ux\stackrel{!}{=}0, \qquad 1\leq \ell\leq n.
\end{equation}
\bl
For all $(\widetilde{w}_0,\widetilde{w}_h,a_0,\veps)\in X_{\gamma-2}\times\R^2 $ sufficiently small, there exists a unique solution to $h_m=0,$ $0\leq m\leq n$, given through the smooth functions
\begin{equation}\label{e:aexp0}
    \alpha=r_\alpha(\widetilde{w}_0,\widetilde{w}_h,a_0,\veps),\qquad a_j=r_j(\widetilde{w}_0,\widetilde{w}_h,a_0,\veps),\ 1\leq j\leq n.
\end{equation}
Moreover, $r_\alpha$,$(r_j)_{1\leq j\leq n}$, and their derivatives with respect to $\widetilde{w}_0,\widetilde{w}_h$, and $a_0$  vanish at the origin. The derivatives of the corrections $a_j$ with respect to $\veps$ at the origin are given through
\begin{equation}\label{e:aexp}
\partial_\veps a_j(0)=\frac{2\pi}{(d_{\|}d_\perp^{n-1})^{1/2}\|u_*^\prime\|_{L^2(\T_{2\pi})}^2}
\int_{\R^n}H_{1,j}\caG\rmd\ux,\qquad 1\leq j\leq n.
\end{equation}
\el
\begin{Proof}
First, notice that all integrals define smooth functions in 
\[
(\widetilde{w}_0,\widetilde{w}_h,\ua,\alpha,a_0,\veps)\in X_{\gamma-2}\times \R^{n+3},
\]
in a small neighborhood of the origin. In order to solve for $(\ua,\alpha)$, it is therefore sufficient to show invertibility of the Jacobi matrix, $A=(A_{\ell j})_{(n+1)\times(n+1)}$, 
 with entries $\partial_\alpha h_j$, $\partial_{a_\ell}h_j$ evaluated at the origin.
    %
More specifically, we note that 
\[
\begin{aligned}
\int_{\R^n}H_0f\rmd\ux&=\langle \caB f, e\rangle\bigg|_{\unu=0}=\int_{\R^n}H_0P_{0,2r}f\rmd\ux, \\
\int_{\R^n}H_{1,\ell}f\rmd\ux&=\partial_{\rmi\nu_\ell}\bigg(\langle \caB f, e\rangle\bigg)\bigg|_{\unu=0}=\int_{\R^n}H_{1,\ell}P_{0,2r}f\rmd\ux, \quad 1\leq \ell\leq n,
\end{aligned}
\]
and the Jacobi matrix has entries
 \[
 \begin{aligned}
 \ell=j=0: \qquad & A_{00}=\int_{\R^n}H_0P_{0,2r}\phi\rmd \ux=\int_{\R^n} H_0 \phi\rmd \ux=\int_{\R^n} \chi_R(\ux) \Big(u^\prime(x_1)\Big)^2\rmd \ux>0;\\
 \ell=0, 1\leq j\leq n: \qquad& A_{0 j}:=\int_{\R^n} H_{0}P_{0,2r}  \partial_{a_j} \caR(\psi)\mid_{\ua=0}\rmd \ux=\int_{\R^n}H_0 \caL_*\left(u_*^\prime\psi_j+u_{*,k}(\psi_j)_{x_1}\right)\rmd \ux;\\
     1\leq \ell\leq n, j=0: \qquad & A_{\ell 0}=\int_{\R^n} H_{1,\ell} P_{0, 2r}\phi\rmd \ux=\int_{\R^n} H_{1,\ell} \phi\rmd \ux;\\
 1\leq \ell, j\leq n: \qquad&
    A_{\ell j}:=\int_{\R^n} H_{1,\ell} P_{0,2r} \partial_{a_j} \caR(\psi)\mid_{\ua=0}\rmd \ux=\int_{\R^n}H_{1,\ell} \caL_*\left(u_*^\prime\psi_j+u_{*,k}(\psi_j)_{x_1}\right)\rmd \ux,
 \end{aligned}
 \]
 where we used that $\phi=\chi_R(\ux)u_*'(x_1)$ for some $R>0$ as defined in \eqref{e:phi}.  
 In order to show that $A$ is invertible, we show that $A$ is in fact diagonal and calculate the nonzero diagonal entries. 
 
 \paragraph{Off-diagonal elements: $A_{\ell j}= 0$ for  $\ell\neq j$.} 
This fact follows readily from parity considerations as follows:
 \begin{itemize}
    \item  $\ell=0,j>0$: the integrand is odd in $x_j$ so that the integral vanishes;
    \item  $\ell=1$: in both cases $j>1$ and $j=0$, the integrand is again odd in $x_1$ and the integral vanishes;
    \item  $\ell>1$: the integrand is odd in $x_\ell$ and the integral vanishes.
    \end{itemize}
    
\paragraph{Diagonal elements:  $A_{\ell \ell}\neq  0$ for all $\ell$.}
We already noted that $A_{00}>0$. For $\ell\geq 1$, we'd like to integrate by parts, which yields an integral $\int_{\R^n}(\caL_*H_{1,\ell})(u_*^\prime\psi_\ell+u_{*,k}(\psi_\ell)_{x_1})\rmd\ux,$ that of course vanishes since $\caL_*H_{1,\ell} =0$. The coefficients $A_{\ell\ell}$ are therefore given by boundary terms in the partial integration. To compute those boundary terms, we repeatedly use the divergence theorem,
\[
\begin{aligned}
   & \int_\Omega (f\Delta g-g\Delta f )\rmd\mathbf{x} =\int_{\partial \Omega}(f\nabla g-g\nabla f)\cdot\mathbf{n}\rmd S,\\
   &\int_\Omega (f\Delta^2 g-g\Delta^2 f )\rmd\mathbf{x} =\int_{\partial \Omega}\Big[f\nabla (\Delta g)-(\Delta g)\nabla f+(\Delta f)\nabla g-g\nabla(\Delta f)\Big]\cdot\mathbf{n}\rmd S,
\end{aligned}
\]
to find that, for any $f\in \ker(\caL_*)$,
\[
\int_\Omega f\caL_* g\rmd \mathbf{x}=-\int_{\partial \Omega}\Big[f\nabla (\Delta g)-(\Delta g)\nabla f+(\Delta f)\nabla g-g\nabla(\Delta f)+2\Big(f\nabla g-g\nabla f\Big)\Big]\cdot\mathbf{n}\rmd S.
\]
Since the integrands are in $L^1$, we can find the integral on $\R^n$ by passing to the limit $R\to\infty$ in domains  $\Omega_R=[-R,R]^n$. Denote the normalized vector in the $x_j$ direction as $\mathbf{n}_j$, and define
\[
(\partial\Omega_R)_{j,\pm}=\left\{ \ux\in\Omega_R\,\middle|\, x_j=\pm R \right\},
\]
so that $\partial\Omega_R=\bigcup_{1\leq j\leq n}\bigcup_{k=\pm}(\partial\Omega_R)_{j,k}$. As a result, we have 
{\small\[
    \int_{\Omega_R} f\caL_* g\rmd \mathbf{x}=-\sum_{\substack{1\leq j\leq n \\ k=\pm}}k
\int_{(\partial \Omega_R)_{j,k}}\overbrace{\Big[f(\Delta g)_{x_j}-(\Delta g) f_{x_j}+(\Delta f)g_{x_j}-g(\Delta f)_{x_j}+2\Big(f g_{x_j}-g f_{x_j}\Big)\Big]}^{:=w_j(\ux)}\rmd S.
\]}
Altogether, we find that for  $1\leq \ell\leq n$, we have 
\begin{equation}\label{e:All}
    A_{\ell\ell}=\lim_{R\to\infty}\int_{[-R,R]^n}f_\ell\caL_* g_\ell\rmd \ux =\lim_{R\to\infty}\left(-\sum_{\substack{1\leq j\leq n \\ k=\pm}}k
\int_{(\partial \Omega_R)_{j,k}} w_{\ell,j}(\ux)\rmd S\right),
\end{equation}
where $f_\ell=H_{1,\ell}$, $g_\ell=u_*^\prime\psi_\ell+u_{*,k}(\psi_\ell)_{x_1}$, and 
\[
w_{\ell,j}:=f_\ell(\Delta g_\ell)_{x_j}-(\Delta g_\ell) (f_\ell)_{x_j}+(\Delta f_\ell)(g_\ell)_{x_j}-g_\ell(\Delta f_\ell)_{x_j}+2\Big(f_\ell (g_\ell)_{x_j}-g_\ell (f_\ell)_{x_j}\Big).
\]
The parities of $f_\ell$ and $g_\ell$ with respect to coordinates $x_j$ simplify the expression of $A_{\ell\ell}$. Specifically, for $\ell=1$, $f_1$ and $g_1$ are even in all $x_k$ directions while, for $\ell>1$, $f_\ell$ and $g_\ell$ are odd in $x_1$ and $x_\ell$, and even in the remaining $x_k$ directions. Both cases lead to the fact that $w_{\ell, j}$ is odd in the $x_j$ direction and even in the remaining $x_k$ directions, which leads to the simplified form
    \begin{equation}\label{e:All-red}
        \int_{[-R,R]^n}f_\ell\caL_*g_\ell\rmd\ux=-2\sum_{1\leq j\leq n}\int_{(\partial\Omega_R)_{j,+}} w_{\ell,j}(\ux)\rmd S.
    \end{equation}
Noting that the terms of order $\caO(|\ux|^{-n})$ or higher in the expression for $w_{\ell, j}$ converge to zero in the estimate of $A_{\ell\ell}$ as $R$ goes to infinity, we derive leading order estimates of $w_{\ell,j}$ as $|\ux|$ goes to infinity. More specifically, for sufficiently large $|\ux|$, we denote $G_\ell:=\partial_{x_\ell}G$ and have 
\[
\begin{aligned}
    g_\ell&=u_*^\prime G_\ell+u_{*,k}G_{1\ell},\\
    (g_\ell)_{x_j}&=\begin{cases}
        u_*^{\prime\prime}G_\ell+(u_*^\prime+u_{*,k}^\prime)G_{1\ell}+\caO(|x|^{-n-1}), & j=1;\\
        u_*^\prime G_{\ell j}+\caO(|\ux|^{-n-1}), & j>1,
    \end{cases}\\
    \Delta g_\ell&=u_*^\tprime G_\ell+(2u_*^\dprime+u_{*,k}^\dprime)G_{1\ell}+\caO(|\ux|^{-n-1}),\\
    (\Delta g_\ell)_{x_j}&=\begin{cases}
    u_*^{(4)} G_\ell+(3u_*^\tprime+u_{*,k}^\tprime)G_{1\ell}+\caO(|\ux|^{-n-1}), & j=1;\\
    u_*^\tprime G_{\ell j}+\caO(|\ux|^{-n-1}), & j>1.
    \end{cases}
\end{aligned}
\]
Now, for $j=1$, we find
    \[
    \begin{aligned}
        w_{\ell,1}=&
        \left[ \sum_{m=0}^3(-1)^m(\partial_{x_1}^{m}f_\ell)u_*^{(4-m)}+2\sum_{\tau=0}^1(-1)^\tau (\partial_{x_1}^{\tau}f_\ell)u_*^{(2-\tau)}\right]G_\ell\\
        &+\left\{ \sum_{m=0}^3(-1)^m(\partial_{x_1}^{m}f_\ell)\left[(3-m)u_*^{(3-m)}+u_{*,k}^{(3-m)}\right]\right.\\
    &+\left.2\sum_{\tau=0}^1(-1)^\tau (\partial_{x_1}^{\tau}f_\ell)\left[(1-\tau)u_*^{(1-\tau)}+u_{*,k}^{(1-\tau)}\right]\right\}G_{1\ell}+\caO(|\ux|^{-n})\\
    =&\overbrace{\left[ 4u_*^\prime u_*^\tprime-2(u_*^\dprime)^2+2(u_*^\prime)^2-\sum_{m=1}^4(-1)^{m} u_*^{(m)}u_{*,k}^{(4-m)}-2\sum_{\tau=1}^2(-1)^\tau u_*^{(\tau)}u_{*,k}^{(2-\tau)}\right]}^{:=W_{\|}(x_1)}\\
            & \cdot\left(x_\ell G_{\ell 1}-\delta_{\ell 1}G_1 \right)+\caO(|\ux|^{-n}).
    \end{aligned}
    \]
    Clearly,  $W_{\|}(x_1)$ is $2\pi$-periodic. We claim that $W_{\|}(x_1)$ is in fact constant. For this,  a direct calculation shows that
    \begin{equation}
    \label{e:Wprime}
         W_{\|}^\prime(x_1)=4u_*^\prime(u_*^{(4)}+u_*^\dprime)-u_*^{(5)}u_{*,k}+u_*^\prime u_{*,k}^{(4)}-2(u_*^\tprime u_{*,k}-u_*^\prime u_{*,k}^\dprime).
    \end{equation}
    On the other hand, $u_*$ solves the $4$-th order ODE
    \[
    -(1+k^2\partial_{x_1}^2)^2u+\mu u-u^3=0,
    \]
    at $k=1$. Taking  derivatives of the above ODE with respect to $x_1$ and $k$ respectively at $u=u_*$ and $k=1$, yields
    \begin{subequations}
        \begin{align}
        -u_*^{(5)}-2u_*^{\tprime}-(1-\mu+3u_*^2)u_*^\prime&=0, \label{e:uprime}\\
        -u_{*,k}^{(4)}-2u_{*,k}^\dprime-(1-\mu+3u_*^2)u_{*,k}&=4(u_*^{(4)}+u^{\dprime}).\label{e:uk}
    \end{align}
    \end{subequations}
    We multiply \eqref{e:uprime} by $u_{*,k}$, \eqref{e:uk} by $u_*$ respectively, and take the difference to find
    \begin{equation}
    \label{e:uk-uprime}
        u_{*,k}(u_*^{(5)}+2u_*^{\tprime})-u_*(u_{*,k}^{(4)}+2u_{*,k}^\dprime)=4u^\prime(u_*^{(4)}+u^{\dprime}).
    \end{equation}
    Comparing \eqref{e:Wprime} with \eqref{e:uk-uprime}, we can conclude that
    $
    W_{\|}^\prime(x_1)\equiv 0,
    $
    and $W_{\|}(x_1)$ is constant with 
    \begin{equation}
        \begin{aligned}
            W_{\|}(x_1)\equiv & \frac{1}{2\pi}\int_0^{2\pi}W_{\|}(x_1)\rmd x_1 
            =\frac{1}{2\pi} \int_0^{2\pi} \left[-6 (u_*^\dprime)^2+2(u_*^\prime)^2-4u_*^\dprime u_{*,k}^\dprime+ 4u_*^\prime u_{*,k}^\prime\right]\rmd \xi=-\frac{\|u_*^\prime\|_{L^2(\T_{2\pi})}^2}{2\pi}d_{\|}.
        \end{aligned}
    \end{equation}
    which simplifies the expression of $w_{\ell, 1}$; that is,
    \[
    w_{\ell, 1}=-\frac{\|u_*^\prime\|_{L^2(\T_{2\pi})}^2}{2\pi}d_{\|} \left(x_\ell G_{\ell 1}-\delta_{\ell 1}G_1 \right)+\caO(|\ux|^{-n}).
    \]
    Noting that $G_\ell(R\ux)=R^{n-1}G_\ell(\ux)$ and $G_{\ell j}(R\ux)=R^{-n}G_{\ell j}(\ux)$, we conclude that
    \begin{equation}\label{e:wl1}
    \begin{aligned}
        \lim_{R\to\infty} \int_{(\partial\Omega_R)_{1,+}}w_{\ell,1}\rmd S=&-\frac{\|u_*^\prime\|_{L^2(\T_{2\pi})}^2}{2\pi}d_{\|} \lim_{R\to\infty}\int_{(\partial\Omega_R)_{1,+}}\left(x_\ell G_{\ell 1}-\delta_{\ell 1}G_1 \right)\rmd S\\
        =&-\frac{\|u_*^\prime\|_{L^2(\T_{2\pi})}^2}{2\pi}d_{\|}\int_{(\partial\Omega_1)_{1,+}}\left(x_\ell G_{\ell 1}-\delta_{\ell 1}G_1 \right)\rmd S.
    \end{aligned}
    \end{equation}
 Next, for $j>1$,
    \[
     \begin{aligned}
     w_{\ell,j}=&\left[f_\ell u_*^\tprime+(\partial_{x_1}^2f_\ell) u_*^\prime+2f_\ell u_*^\prime \right]G_{\ell j}-\delta_{\ell j}2u_*^\prime(u_*^\tprime+u_*^\prime)G_\ell
     +\caO(|\ux|^{-n})\\
     =&
        \overbrace{2[u_*^\prime u_*^\tprime+(u_*^\prime)^2]}^{:=W_\perp(x_1)}(x_\ell G_{\ell j}-\delta_{\ell j}G_j)+\caO(|\ux|^{-n}). 
     \end{aligned}
     \]
     We note that $W_\perp$ is $2\pi$-periodic in $x_1$ and calculate its average
     \[
     \begin{aligned}
         \overline{W_{\perp}}=&\frac{1}{2\pi}\int_0^{2\pi}2[u_*^\prime u_*^\tprime+(u_*^\prime)^2]\rmd x_1=\frac{1}{2\pi}\int_{\T_{2\pi}}2[-(u_*^\dprime)^2+(u_*^\prime)^2]\rmd \xi=-\frac{\|u_*^\prime\|_{L^2(\T_{2\pi})}^2}{2\pi}d_{\perp}.
     \end{aligned}
     \]
     Introducing the $2\pi$-periodic function 
     \[
     W_0(x_1):=\int_0^{x_1}(W_\perp(\xi)-\overline{W_\perp})\rmd \xi,
     \]
     and $(\partial\Omega_R)_{j,+}^\perp:=\{\ux_\perp\,\mid\, (x_1,\ux_\perp)\in (\partial\Omega_R)_{j,+}\}$,
     we have $W_\perp=\overline{W_\perp}+W_0^\prime$ and
     \begin{equation}\label{e:wlj}
         \begin{aligned}
             &\lim_{R\to\infty} \int_{(\partial\Omega_R)_{j,+}}w_{\ell,j}\rmd S\\
             &\qquad= \overline{W_\perp}\lim_{R\to\infty}\int_{(\partial\Omega_R)_{j,+}}(x_\ell G_{\ell j}-\delta_{\ell j}G_j)\rmd S+
             \lim_{R\to\infty}\int_{(\partial\Omega_R)_{j,+}} W_0^\prime (x_\ell G_{\ell j}-\delta_{\ell j}G_j)\rmd S\\
             &\qquad=\overline{W_\perp}\int_{(\partial\Omega_1)_{j,+}}(x_\ell G_{\ell j}-\delta_{\ell j}G_j)\rmd S\\
              &\qquad\qquad+
             \lim_{R\to\infty}\left[W_0\overbrace{\int_{(\partial\Omega_R)_{j,+}^\perp}(x_\ell G_{\ell j}-\delta_{\ell j}G_j)\rmd S\Bigg|_{x_1=-R}^{x_1=R}}^{=\caO(1/R)} -\overbrace{\int_{(\partial\Omega_R)_{j,+}} W_0\partial_{x_1} (x_\ell G_{\ell j}-\delta_{\ell j}G_j)\rmd S}^{=\caO(1/R)}\right]\\
             &\qquad =-\frac{\|u_*^\prime\|_{L^2(\T_{2\pi})}^2}{2\pi}d_{\perp}\int_{(\partial\Omega_1)_{j,+}}(x_\ell G_{\ell j}-\delta_{\ell j}G_j)\rmd S.
         \end{aligned}
     \end{equation}
Recall that  the Green's function $G$ is a function of the anisotropic radial variable $r=\sqrt{\frac{x_1^2}{d_{\|}}+\frac{|\ux_\perp|^2}{d_\perp}}$ and that $G^\prime:=\frac{\rmd G}{\rmd r}=-\frac{r^{1-n}}{C_n}$ where $C_n>0$ is the surface area of the unit ball in $\R^n$.  Slightly simplifying notation, we write $d_1:=d_{\|}$, $d_\ell:=d_\perp$ for $\ell>1$, and have
    \begin{equation}\label{e:G-int}
        \begin{aligned}
 x_\ell G_{\ell j}-\delta_{\ell j}G_j=&x_\ell\left[\frac{1}{d_\ell d_j}\frac{x_\ell x_j}{r^2}G^\dprime+\left(\frac{\delta_{\ell j}}{d_j}\frac{1}{r}-\frac{1}{d_\ell d_j}\frac{x_\ell x_j}{r^3}\right)G^\prime\right]-\delta_{\ell j}\frac{1}{d_j}\frac{x_j}{r}G^\prime\\
 =&\frac{1}{d_\ell d_j}\frac{x_\ell^2 x_j}{r^3}(r G^\dprime- G^\prime)=\frac{n C_n}{d_\ell d_j}\frac{x_\ell^2 x_j}{r^{n+2}}.
    \end{aligned}
    \end{equation}
Combining \eqref{e:All}, \eqref{e:All-red},\eqref{e:wl1}, \eqref{e:wlj}, and \eqref{e:G-int}), we have
\begin{equation}
A_{\ell\ell}=\frac{n\|u_*^\prime\|_{L^2(\T_{2\pi})}^2}{\pi C_n d_\ell}\left(\sum_{j=1}^n\int_{(\partial\Omega_1)_{j,+}}\frac{x_\ell^2 }{r^{n+2}}\rmd S\right)>0.
\end{equation}
More specifically, defining 
\[
\widetilde{x}_1:=\frac{x_1}{d^{1/2}_{\|}}, \quad
\widetilde{\ux}_\perp:=\frac{\ux_\perp}{d^{1/2}_{\perp}}, \quad 
\widetilde{r}:=|\widetilde{x}|,\quad
\widetilde{\Omega}_1:=\{\widetilde{x}\mid (d^{1/2}_{\|}\widetilde{x}_1,d^{1/2}_{\perp}\widetilde{\ux}_\perp)\in\Omega_1\},
\]
we find by rescaling,  using the divergence theorem, and replacing the volume of integration by the unit ball,   
\begin{equation}
\begin{aligned}
A_{\ell\ell}=&\frac{n\|u_*^\prime\|_{L^2(\T_{2\pi})}^2}{\pi C_n d_\ell}\left(\sum_{j=1}^n\int_{(\partial\Omega_1)_{j,+}}\frac{x_\ell^2 }{r^{n+2}}\rmd S\right)\\
=&\frac{n\|u_*^\prime\|_{L^2(\T_{2\pi})}^2}{\pi C_n}\left[d_{\perp}^{(n-1)/2}\int_{(\partial\widetilde{\Omega}_1)_{1,+}}\frac{\widetilde{x}_\ell^2 }{\widetilde{r}^{n+2}}\rmd S+(d_{\|}d_\perp^{n-2})^{1/2}\sum_{j=2}^n\int_{(\partial\widetilde{\Omega}_1)_{j,+}}\frac{\widetilde{x}_\ell^2 }{\widetilde{r}^{n+2}}\rmd S\right]\\
=& \frac{n\|u_*^\prime\|_{L^2(\T_{2\pi})}^2}{2\pi C_n}(d_{\|}d_\perp^{n-1})^{1/2}\int_{\partial\widetilde{\Omega}_1}\frac{\widetilde{x}_\ell^2 }{\widetilde{r}^{n+2}}\widetilde{\ux}\cdot \mathbf{n}\rmd S\\
=&\frac{n\|u_*^\prime\|_{L^2(\T_{2\pi})}^2}{2\pi C_n}(d_{\|}d_\perp^{n-1})^{1/2}\left[\int_{\widetilde{\Omega}_1\backslash B_r(0)}\nabla\cdot\left(\frac{\widetilde{x}_\ell^2 }{\widetilde{r}^{n+2}}\widetilde{\ux}\right)\rmd V+\int_{\partial B_r(0)}\left(\frac{\widetilde{x}_\ell^2 }{\widetilde{r}^{n+2}}\widetilde{\ux}\right)\cdot \mathbf{n}\rmd S
\right]\\
=&\frac{n\|u_*^\prime\|_{L^2(\T_{2\pi})}^2}{2\pi C_n}(d_{\|}d_\perp^{n-1})^{1/2}\int_{\partial B_1(0)}x_\ell^2\rmd S\\
=&\frac{\|u_*^\prime\|_{L^2(\T_{2\pi})}^2}{2\pi C_n}(d_{\|}d_\perp^{n-1})^{1/2}\int_{\partial B_1(0)}|\ux|^2\rmd S\\
=&\frac{\|u_*^\prime\|_{L^2(\T_{2\pi})}^2}{2\pi }(d_{\|}d_\perp^{n-1})^{1/2},
\end{aligned}
\end{equation}
which concludes the proof.
\end{Proof}

\subsection{Solving the equation on the range and finding the pinning location --- conclusion of proof}
We are now ready to solve \eqref{e:ift-form-ref2}, that is,
\[
0=-\begin{pmatrix}
\widetilde{w}_0 \\ \widetilde{w}_h
\end{pmatrix}+\begin{pmatrix}
\widetilde{P}_{0,r} \\
\widetilde{P}_{h,r}
 \end{pmatrix} \caL_*^{-1}\left(\caR(\psi(\underline{a})) + \caN(\psi(\underline{a}), \widetilde{w}_0 u_*^\prime +\widetilde{w}_h)+\caG(a_0, \veps)+\alpha \phi\right),
\] 
 after substituting the functions $r_\alpha$ and $r_j$,  for $\alpha$ and $a_j$, respectively, $1\leq j\leq n$. The linear operator $\caL_*$ is bounded invertible on its range and the linearization of the right-hand side equals to the negative identity.

It remains to find solutions where the dummy correction $\alpha$ vanishes. 
The parameter $\alpha$ is found readily by projection of the equation onto $H_0=u_*'$. 
We find at leading order 
\begin{equation}\label{e:x1*}
\alpha(\veps;{a_0})=\veps \frac{1}{A_{00}}\int_{\R^n} g(\ux)u_*'(x_1+{a_0})\rmd \ux +\rmO(\veps^2)=\veps \left(\frac{M({a_0})}{A_{00}}+\rmO(\veps)\right),
\end{equation}
where $M$ was defined in \eqref{e:melnikov}. We can solve $\alpha(\veps;{a_0})=0$ for ${a_0}={a_0}^*(\veps)$ near {${a_0}=x_1^*$}, defined in Hypothesis \ref{h:2}, using that $M'(x_1^*)\neq 0$ as assumed in Hypothesis \ref{h:2}. We therefore solve  after dividing by $\veps$ with implicit function theorem. Similarly, we find $\ua$ from \eqref{e:aexp} after substituting the value of $x_1^*$ found  from \eqref{e:x1*}.
This concludes the proof of Theorem~\ref{t:main}. We note here that the scheme is slightly simpler for positive $\veps$ since we can use $a_0=x_1^*$ as a variable rather than the dummy variable $\alpha$.
\begin{corollary}\label{c:inv}
For all $\veps$ sufficiently small, the linearization of 
\begin{align}
0&=-\begin{pmatrix}
\widetilde{w}_0 \\ \widetilde{w}_h
\end{pmatrix}+ \caL_*^{-1}\begin{pmatrix}   
\widetilde{P}_{0,r} (1-Q)\\
\widetilde{P}_{h,r}
 \end{pmatrix}\left(\caR(\psi(\underline{a})) + \caN(\psi(\underline{a}), \widetilde{w}_0 u_*^\prime +\widetilde{w}_h)+\caG(a_0, \veps)\right),\label{e:ift0}\\
 0&=\mathcal{Q}\left(\caR(\psi(\underline{a})) + \caN(\psi(\underline{a}), \widetilde{w}_0 u_*^\prime +\widetilde{w}_h)+\caG(a_0, \veps)\right),\label{e:ift}
\end{align}
with respect to the variables $(\widetilde{w}_0, \widetilde{w}_h,\ua,a_0)$ is bounded invertible on $R_{\gamma-1}\times R_\gamma\times \R^n\times\R$. Here, $\widetilde{P}_{0,r}$ and $\widetilde{P}_{h,r}$ denote the mode filters and $\mathcal{Q}$ denotes the projection onto the cokernel of $\caL_*$ introduced in \eqref{e:caQ}, that is, $\mathcal{Q}w=0$ when the scalar products with the pseudo-harmonic polynomials of degree 0 and 1 vanish; compare \eqref{e:caQ}. 
\end{corollary}
Note that the range of $\mathcal{Q}$ is in fact irrelevant, as \eqref{e:ift} is independent of choosing a particular basis in the complement of the range, while solving \eqref{e:ift} makes the appearance of $\mathcal{Q}$ in \eqref{e:ift0} irrelevant. 

\section{Discussion}\label{s:6}

We studied the effect of localized impurities on striped phases. In order to capture the effect of essential spectrum in the linearization which prevents a straightforward use of the implicit function theorem, we developed a description of the linearization via a conjugation to a ``homogenized'' linear problem, that turns out to be an anisotropic Laplacian, together with an invertible part acting on finite-wavelength perturbations. The second ingredient was an explicit phase modulation that compensates for the cokernel of the Laplacian. The resulting equation is then well-posed on the natural algebraically weighted spaces and can be solved with the Implicit Function Theorem. 

The analysis here points to a wealth of interesting extensions. First, setting the problem in the context of an Implicit Function Theorem naturally suggests continuation in parameter space, for not necessarily small $\veps$, up until possible bifurcations. It also provides a straightforward algorithm toward numerical approximations. In this latter direction, it would be interesting to obtain  higher-order, multi-pole farfield corrections theoretically and computationally, improving on the $\caO(|\ux|^{-(n-1)})$-approximation through the modulations $\psi$ used here. 

In a different direction, one notices that the perturbation here is not critical in a scaling sense, that is, the leading order corrections are in fact determined by the linear part. It would be interesting to investigate if this is more generally true. In this spirit, a particularly relevant case arises when the wavenumber of the stripes is zigzag-critical. These wavenumbers turn out to be energy-minimizing among periodic patterns and are selected by grain boundaries and possibly other defects \cite{LloydScheel}. The linearization  at these energy-minimizing stripes is marginally stable, with a Fourier-Bloch long-wavelength expansion of the form $\lambda\sim \nu_x^2-\nu_y^4$ which gives slower decay of the Green's function; see \cite{hile2,hile} for Fredholm results that mimic Theorem \ref{t:McOwen} in that situation. The decay is sufficient to guarantee temporal decay when the inhomogeneity is localized in space and time \cite{chowdhary2024weak} and in fact subcritical in that it does not alter linear asymptotics. It would be interesting to establish the corresponding subcriticality also in the present case of a time-stationary inhomogeneity. A related direction that would probe criticality are inhomogeneities with weaker localization, such as for instance $g(x_1)g(x_2)$ with $g(x_2)$ bounded and $g(x_1)$ localized. Again, one would hope to establish that the linear problem predicts correct farfield asymptotics.  

Finally, one may wish to reintroduce time-dependence. Periodic time-dependence in $g$ should lead to rather equivalent results, with only the time-average of $g$ contributing. More interestingly, one would like to understand stability of profiles studied here. One may suspect that the sign of $M^\prime$ indicates stability and, in the simplest case when $M$ has two nondegenerate zeros, one of the associated pinned solutions is unstable  with a positive eigenvalue in the linearization while the other zero gives rise to a stable solution. Such pinning dynamics are well understood in the absence of continuous spectrum, for instance when analyzing front pinning in heterogeneous media \cite{st}. They are also readily analyzed via center manifold reduction when inhomogeneities are coperiodic with the pattern. In the present case, adjusting the position of the periodic pattern with respect to the inhomogeneity is possible only locally uniformly, and one would expect that this adjustment gives rise to a diffusive mixing wave in the far field; see \cite[Fig. 4]{LloydScheel} for a numerical illustration of this selection in the case of grain boundaries and \cite{SSSU_12} for an analysis of the mixing of phases. More specifically, simply proving  stability of the solution without unstable eigenvalues in the linearization would be first step towards understanding these pinning dynamics, that is, how striped patterns eventually converge towards the pinned solutions found here.

\end{document}